\documentstyle[bezier,12pt]{article}
\input{emlines.sty}
\input{amssym.def}  
\input{amssym}      
\newcommand{\mathbb}[1]{{\Bbb #1}}
\newcommand{\mathcal}[1]{{\cal #1}}
\newcommand{\mathrm}[1]{\mbox{#1}}
\newcommand{\emph}[1]{{\em #1\/}}

\newcommand{\ensuremath}[1]{#1}
\newcommand{\textit}[1]{{\it #1}}

\newif\ifnotesw\noteswtrue
   {\ifnotesw\marginpar[\hfill\(\top\)]{\(\top\)}\fi}%
   {\ifnotesw\marginpar[\hfill\(\bot\)]{\(\bot\)}\fi}

\newcommand{\mnote}[1]%
    {\ifnotesw\marginpar%
        [{\scriptsize\it\begin{minipage}[t]{\marginparwidth}\raggedleft#1%
                        \end{minipage}}]%
        {\scriptsize\it\begin{minipage}[t]{\marginparwidth}\raggedright#1%
                        \end{minipage}}%
    \fi}
\noteswfalse   

\newcommand{\qed}{$\;\;\;\Box$}
\newenvironment{proof}{\par\smallbreak{\sl Proof.~}}
{\unskip\nobreak\hfill \qed \par\medbreak}

\newcommand{\hide}[1]{}

\renewcommand{\L}{\ensuremath{\mathrm{L}}}
\newcommand{\Con}{\ensuremath{\mathrm{C}}}

\newcommand{\D}{\ensuremath{{\cal D}}}

\newcommand{\mb}[1]{\ensuremath{\mathbb{#1}}}

\newcommand{\R}{\mb{R}}

\newcommand{\G}{\ensuremath{{\cal G}}}

\newcommand{\mes}{\mathrm{mes}}

\renewcommand{\d}{\ensuremath{\partial}}

\newfont{\bl}{msbm10 scaled \magstep2}

\newtheorem{thm}{Theorem}

\newtheorem{lemma}{Lemma}
\newtheorem{prop}{Proposition}
\newtheorem{defn}{Definition}
\newtheorem{cor}{Corollary}

\newtheorem{rem}{Remark}

\newtheorem{claim}{Claim}
\newtheorem{cla}{Claim}
\newtheorem{clai}{Claim}
\newtheorem{fact}{Fact}

\newcommand{\bem}[1]{\vadjust{\rlap{\kern\hsize\thinspace\vbox%
                       to0pt{\hbox{$\clubsuit${\small\tt #1}}\vss}}}}
\newcommand{\add}{\vadjust{\rlap{\kern\hsize\thinspace\vbox%
                       to0pt{\hbox{\enspace\enspace$\clubsuit$}\vss}}}}

\newcommand{\notmid}{\mid\kern-0.5em\not\kern0.5em}

\newcommand{\de}{\delta}
\newcommand{\eps}{\varepsilon}

\newcommand{\vphi}{\varphi}

\newcommand{\la}{\lambda}
\newcommand{\om}{\omega}

\newcommand{\supp}{\mathop{\mathrm{supp}}}

\newcommand{\sing}{\mathop{\mathrm{sing}}}

\renewcommand{\em}[1]{{\it #1\/}}

\title{Delta Waves
for a Strongly Singular Initial-Boundary
Hyperbolic Problem with Integral Boundary Condition
}

\newcounter{thesame}
\author{
I.~Kmit\thanks{
Department of Numerical Mathematics \& Programming,
National University ``Lvivska Polytechnika'',
Bandera St.\ 12, 290646 Lviv, Ukraine. E-mail:
{\tt kmit@ov.litech.net}
This work was done
while visiting the Institut f{\"u}r Mathematik,
Universit{\"a}t Wien,
supported by an {\"O}AD grant.}
}

\date{22 July 2003}

\begin{document}

\maketitle

\begin{abstract}
We investigate the existence and the singular structure
of delta wave solutions to a semilinear strictly hyperbolic equation with strongly singular
initial and boundary conditions.
The boundary conditions are given in nonlocal form with a linear integral operator involved.
We construct a delta wave solution as a distributional limit
of solutions to the
regularized system.
This determines  the macroscopic behavior of the corresponding
generalized solution in the Colombeau algebra
$\G$ of generalized functions.
We represent our delta wave as a sum of a purely singular part
satisfying a linear system and a regular part
satisfying a nonlinear system.


\end{abstract}

\section{Introduction}\label{sec:intr}

In the domain
$$
\Pi=\{(x,t)\in\R^2\,|\,0<x<L,t>0\}
$$
we study the following initial-boundary
value problem for the first-order semilinear hyperbolic equation:
\begin{eqnarray}
(\partial_t  + \lambda(x,t)\partial_x) u &=& p(x,t)u + f(x,t,u), \qquad (x,t)\in \Pi
  \label{eq:1}\\
u|_{t=0} &=& a(x), \qquad x\in (0,L) \label{eq:2} \\
u|_{x=0} &=& \int\limits_{0}^{L}d(x,t)u\,dx,
  \qquad  t\in(0,\infty) \label{eq:3} \, .
\end{eqnarray}
Mathematical models of this kind stem from mathematical biology and serve to describe
the age-dependent population dynamics (see~\cite{bo-ar,bu-ia,Huyer,Sanchez,Sond}).
In particular, the linear case of the
problem, when $f(x,t,u)$ does not depend on $u$,
arises in demography, where  $u(x,t)$
is the population density of age $x$ at time $t$,
$a(x)$ is the initial density,
$d(x,t)$ is the birth rate,
$-p(x,t)$ is the death rate, and $f(x,t)$ is the migrant density.
Nonlinear models of age structured populations are studied
in~\cite{bo-ar,bu-ia}.
To model
point-concentration of the initial density and the birth rate, we consider the data $a(x)$ and $d(x,t)$
to be strongly singular,
of the Dirac delta type.

As well known, solutions to the classical initial-boundary
semilinear hyperbolic problems in a single space variable
are at least as singular as the initial
and the boundary data.
We therefore
 can expect for the nonclassical problem (\ref{eq:1})--(\ref{eq:3})
that the multiplication of distributions  appears in the right-hand
sides of (\ref{eq:1}) and  (\ref{eq:3}).
Such multiplication in general cannot be performed within
the distributional theory and,
by this reason,  is
 usually defined in differential algebras of generalized functions.
In~\cite{KmitHorm} we used the Colombeau algebra of generalized functions
$\G(\overline{\Pi})$  \cite{1,3,11} to prove
a global existence-uniqueness
result for (\ref{eq:1})--(\ref{eq:3}).
Nevertheless, the macroscopic behavior
of the Colombeau solutions remained unclear.

We here show  that  the Colombeau solution to
(\ref{eq:1})--(\ref{eq:3}) is associated
to the distributional solution.
This means that the system has a delta wave
solution in the sense of~\cite{RauchReed87}, i.e.,   the
sequence of approximate
(or sequential) solutions obtained by regularizing  all  singular data
has a weak limit.
In the course of construction of the delta wave solution
we show  interaction and propagation of singularities.

 It should be noted that an associated distribution
or in other terminology, a delta wave solution,
though contains an important information about
the singular structure of the generalized solution,
in general does not satisfy the system in a differential-algebraic sense.
Our paper brings one more example into the collection of associated
distributions which are not distributional solutions.

The advantage of using delta wave solutions lies in the fact that, due to the procedure
of their obtaining, they
 are always  stable.
In contrast with this, if we
use a priori defined intrinsic multiplication of distributions
for
obtaining  distributional solutions, the result may be
 nonstable and noncorrect~\cite{2}.
The concept of a delta wave solution has
also other advantages. It
allows us to solve nonlinear
systems and systems with nonsmooth coefficients, for which
the distributional theory is not well adapted.
For the delta wave solutions of semilinear hyperbolic problems we refer the reader
to the sources
\cite{2,3_1,5} and
\cite{KmitHorm}--\cite{RauchReed87}.


We  split a delta wave into the sum of a regular part satisfying a nonlinear
equation and a singular part satisfying a linear equation. The idea of nonlinear
splitting goes back to \cite{8,RauchReed87,14,3_1}.
An important  feature of the nonlinear splitting suggested here is
  quite strong interdependence of the singular and the regular parts.
A similar phenomenon is discovered in~\cite{KmitHorm} for a nonlocal problem with
nonseparable
boundary conditions, where the singular part of the nonlinear
splitting depends on the regular part.

Delta wave solutions for initial-boundary semilinear hyperbolic problems
were considered in \cite{13,KmitHorm}.
The paper \cite{13} investigates  the existence and structure of delta waves
in a nonlinear boundary value problem
for a second order
hyperbolic equation where the boundary condition is nonlinear
and the nonlinearity is given by
a bounded smooth function.
Both in \cite{13}  and \cite{KmitHorm} the right hand side of
the differential equations is
bounded and in \cite{13} it can be also sublinear with respect to $u$.

%

The paper is organized as follows. In Section 2 we describe in
detail our splitting  of a delta wave solution and state our main result.
The proof is given in Sections~3--8. In particular,
in Section 5 we show that our splitting procedure is correct.
In Section 6 we are concerned with the regular part.
Using the Cauchy criterion of the uniform
convergency, we prove that the family of approximate solutions to the
regular part uniformly converges on any compact subset of
$\overline{\Pi}$. In Section 7 we deal with the singular part
and prove that the sequence of approximate solutions
to the singular part converges in $\D'(\overline{\Pi})$ to a function $v$.
We then show that $v$  actually represents the purely singular part of the
initial problem and is the sum of measures concentrated on characteristic
curves (see Sections 7 and 8).

\section{Interaction and propagation of strong singularities  and
construction of a delta wave solution}
\label{sec:delta}

We first list assumptions that will
be made for  the problem (\ref{eq:1})--(\ref{eq:3}).
\vskip1.0mm
{\it Assumption 1.}
$a(x)=a_s(x)+a_r(x),\quad d(x,t)=b(x)\otimes c(t)=(b_s(x)+b_r(x))\otimes (c_s(t)+c_r(t))$, where
$a_s(x)$, $b_s(x)$, $c_s(t)$ and
$a_r(x)$, $b_r(x)$, $c_r(t)$ are, respectively,
 singular and
regular parts of the functions
 $a(x)$, $b(x)$, and $c(t)$.
\vskip1.0mm
{\it Assumption 2.}
$a_s(x), b_s(x),$ and $c_s(t)$ are the finite sums of the delta functions at points, whose
supports are as follows:\\
$
\supp a_s(x)=\{x_1^*,x_2^*,\dots,x_m^*\},\mbox{where}\,\, 0<x_1^*<\dots<x_m^*<L.\\
\supp b_s(x)=\{x_1,x_2,\dots,x_k\},\mbox{where}\,\, 0<x_1<\dots<x_k<L.\\
\supp c_s(t)=\{t_1,t_2,\dots,t_l\},\mbox{where}\,\, 0<t_1<\dots<t_l.
$
\vskip1.0mm
{\it Assumption 3.}
$b_r(0)=0$ and $b_r(L)=0.$

\vskip1.0mm
{\it Assumption 4.}
$a_r(0)=0$ and $c_r(0)=0.$

\vskip1.0mm
{\it Assumption 5.}
$p,f,a_r,b_r,c_r$ are continuous and $\la$ is continuously differentiable with respect to
all their arguments, $f$ is continuously differentiable with respect to~$u$.
\vskip1.0mm
{\it Assumption 6.}
$\la(x,t)>0 $ for all $(x,t)\in\overline{\Pi}$.


\vskip1.0mm
{\it Assumption 7.}
 $f$ and $\nabla_uf$ are globally bounded with respect to $(x,t)$
varying in compact subsets of $\overline{\Pi}$.
\vskip1.0mm
Assumption 4 serves
to ensure the 0-order compatibility between (\ref{eq:2}) and (\ref{eq:3}).
Note that the assumptions are not restrictive from the viewpoint of applications.

Recall that all characteristics of the differential equation (\ref{eq:1}) are solutions to
the following initial problem for ordinary differential equation:
$$
\frac{d\xi}{d\tau}=\la(\xi(\tau),\tau),\quad \xi(t)=x,
$$
where $(x,t)\in\overline{\Pi}.$
It is well known that, under Assumptions 5 and 6, for every $(x,t)\in\overline{\Pi}$
this problem has a unique $\Con^1$-solution which can be expressed in any of two forms
$\xi=\om(\tau;x,t)$ or $\tau=\tilde\om(\xi;x,t)$.


\begin{figure}
\centerline{
\unitlength=1.00mm
\special{em:linewidth 0.4pt}
\linethickness{0.4pt}
\begin{picture}(110.00,114.00)
\put(25.00,5.00){\makebox(0,0)[cc]{$0$}}
\put(110.00,5.00){\makebox(0,0)[cc]{$L$}}
\put(25.00,114.00){\makebox(0,0)[cc]{$T$}}
\put(30.00,10.00){\framebox(75.00,100.00)[cc]{}}
\emline{57.00}{10.00}{1}{57.00}{110.00}{2}
\put(57.00,5.00){\makebox(0,0)[cc]{$x_1$}}
\bezier{388}(41.00,10.00)(86.00,33.00)(105.00,75.00)
\emline{57.00}{19.00}{3}{54.00}{19.00}{4}
\emline{51.00}{19.00}{5}{48.00}{19.00}{6}
\emline{45.00}{19.00}{7}{42.00}{19.00}{8}
\emline{39.00}{19.00}{9}{36.00}{19.00}{10}
\emline{33.00}{19.00}{11}{30.00}{19.00}{12}
\bezier{472}(30.00,19.00)(78.00,57.00)(99.00,110.00)
\emline{30.00}{62.00}{13}{105.00}{62.00}{14}
\bezier{276}(30.00,62.00)(70.00,89.00)(70.00,110.00)
\emline{57.00}{44.00}{15}{54.00}{44.00}{16}
\emline{54.00}{44.00}{17}{57.00}{44.00}{18}
\emline{57.00}{44.00}{19}{54.00}{44.00}{20}
\emline{51.00}{44.00}{21}{48.00}{44.00}{22}
\emline{45.00}{44.00}{23}{42.00}{44.00}{24}
\emline{39.00}{44.00}{25}{36.00}{44.00}{26}
\emline{36.00}{44.00}{27}{33.00}{44.00}{28}
\emline{33.00}{44.00}{29}{36.00}{44.00}{30}
\emline{33.00}{44.00}{31}{30.00}{44.00}{32}
\bezier{340}(30.00,44.00)(66.00,73.00)(77.00,110.00)
\emline{57.00}{71.00}{33}{54.00}{71.00}{34}
\emline{54.00}{71.00}{35}{57.00}{71.00}{36}
\emline{57.00}{71.00}{37}{54.00}{71.00}{38}
\emline{51.00}{71.00}{39}{47.00}{71.00}{40}
\emline{44.00}{71.00}{41}{41.00}{71.00}{42}
\emline{38.00}{71.00}{43}{35.00}{71.00}{44}
\emline{33.00}{71.00}{45}{30.00}{71.00}{46}
\bezier{180}(30.00,71.00)(46.00,89.00)(48.00,110.00)
\bezier{100}(30.00,34.00)(68.00,68.00)(88.00,110.00)
\put(20.00,19.00){\makebox(0,0)[cc]{$t_1^*=\tilde t_1$}}
\put(20.00,44.00){\makebox(0,0)[cc]{$t_2^*=\tilde t_2$}}
\put(20.00,62.00){\makebox(0,0)[cc]{$t_3^*=t_1$}}
\put(20.00,71.00){\makebox(0,0)[cc]{$t_4^*=\tilde t_4$}}
\put(85.00,38.00){\makebox(0,0)[cc]{$I_+$}}
\put(75.00,71.00){\makebox(0,0)[cc]{$I_+$}}
\put(84.00,98.00){\makebox(0,0)[cc]{$I_-$}}
\put(77.00,103.00){\makebox(0,0)[cc]{$I_+$}}
\put(64.00,99.00){\makebox(0,0)[cc]{$I_+$}}
\put(41.00,94.00){\makebox(0,0)[cc]{$I_+$}}
\put(42.00,10.00){\circle{2.83}}
\put(30.00,62.00){\circle{2.00}}
\put(57.00,62.00){\circle{2.00}}
\put(57.00,62.00){\circle*{2.00}}
\put(57.00,19.00){\circle*{2.00}}
\put(57.00,44.00){\circle*{2.00}}
\put(57.00,71.00){\circle*{2.00}}
\put(57.00,84.00){\circle*{2.00}}
\put(30.00,19.00){\circle{2.00}}
\put(30.00,44.00){\circle{2.00}}
\put(30.00,71.00){\circle{2.00}}
\put(42.00,5.00){\makebox(0,0)[cc]{$x_1^*$}}
\emline{37.00}{44.00}{57}{34.00}{44.00}{58}
\end{picture}
}
\end{figure}


Choose $\eps_0>0$ so small that $x_1^*-\eps_0>0$ and $t_1-\eps_0>0$.
Some additional conditions on $\eps_0$ will be put below.
We will consider
$\eps$ in the range
$0<\eps\le\eps_0$.
\begin{defn}
Let $I_-$ be the union of the characteristics
$\om(t;x_i,t_j)$ for all $i\le k$ and
$j\le l$.
Let $I_-^{\eps}$ be the union of the neighborhoods
$\{(\om(\tau;x_i,t),\tau)\,|\,
\tilde\om(x_i;x_i-\eps, t_j+\eps)
<t<\tilde\om(x_i;x_i+\eps, t_j-\eps)\}$
for all $i\le k$ and
$j\le l$.
\end{defn}





\begin{defn}
Let $I_+=\bigcup_{n\ge 0}I_+[n]$ and, for  $\eps<\eps_0$,
$I_+^{\eps}=\bigcup_{n\ge 0}I_+^{\eps}[n]$,
where
$I_+[n]$ and $I_+^{\eps}[n]$ are
subsets of $\overline{\Pi}$ defined
by induction as follows.

\begin{itemize}
\item
 $I_+[0]$ includes the characteristics $\om(t;x_i^*,0)$
and $\om(t;0,t_j)$
for all $i\le m$ and
$j\le l$
(i.e. $I_+[0]$ is the union of these characteristics).

$I_+^{\eps}[0]$ includes the neighborhoods $\{(\om(\tau;x,0),\tau)\,|\,
x_i^*-\eps<x<x_i^*+\eps\}$ and $\{(\om(\tau;0,t),\tau)\,|\,t_j-\eps<t<t_j+\eps\}$.

\item
Let $n\ge 1$. If $I_+[n-1]$ includes the characteristic
 $\om(t;x_i,\tilde t)$, then
$I_+[n]$ includes  the characteristic $\om(t;0,\tilde t)$.

If $I_+^{\eps}[n-1]$ includes the neighborhood
 $\{(\om(\tau;x_i,t),\tau)\,|\,\tilde t-\eps^-<t<\tilde t+\eps^+\}$,
then
$I_+^{\eps}[n]$ includes  the neighborhood
 $\{(\om(\tau;0,t),\tau)\,|\,\tilde\om(x_i-\eps;x_i,\tilde t-\eps^-)<t
<\tilde\om(x_i+\eps;x_i,\tilde t+\eps^+)\}$.
\end{itemize}
\end{defn}
The set $I_+$ captures the propagation  of all singularities.
For characteristics contributing into
$I_+$
(respectively, $I_+\setminus I_+[0]$),
denote their intersection points   with the axis $x=0$
by $t_1^*,t_2^*,\dots$
(respectively, $\tilde t_{i_1},
\tilde t_{i_2},\dots$).
We assume that
$t_j^*<t_{j+1}^*$ for $j\ge 1$
and
$i_n=p$  for
$\tilde t_{i_n}=t_p^*$.
 Obviously, $\{t_1^*,t_2^*,\dots\}=\{t_1,\dots,t_l\}\cup
\{\tilde t_{i_1},\tilde t_{i_2},\dots\}$.
Let $\eps_i^-(\eps)$ and
$\eps_i^+(\eps)$ be such that $I_+^{\eps}\cap \{(x,t)\in\overline{\Pi}\,|\,x=0\}=
\bigcup_i
\{(0,t)\,|\, t_i^*-\eps_i^-(\eps)
<t<t_i^*+\eps_i^+(\eps)\}$.
If  $t_i^*=t_j$ for some $j\le l$, then $\eps_i^-(\eps)=\eps_i^+(\eps)=\eps$.
Observe that $\lim\limits_{\eps\to 0}\eps_i^-(\eps)=0$ and
$\lim\limits_{\eps\to 0}\eps_i^+(\eps)=0$.

\vskip1.0mm
{\it Assumption 8.}
$
\tilde\om(0;x_i,t_j)\ne t_s^*$, $\om(0;x_i,t_j)\ne x_q^*$
 for all $i\le k$,
$j\le l$, $q\le m$
and $t_s^*<t_l$.
 \vskip1.0mm
This assumption means that no three different singularities caused by the initial
 and the boundary data
hit at the same point. In other words,
neither  points
$(x_q^*,0)$ and $(x_i,t_j)$ nor points $(0,t_s^*)$ and $(x_i,t_j)$
are connected by any
of characteristic curves.
 As a consequence,
there exists $\eps_0$ such that, for each
$\eps\le\eps_0$,
$I_-^{\eps}\cap I_+^{\eps}=\emptyset$.
Assume that $(0,0)\not\in I_-$.
We choose $\eps_0$ so small
that $I_-^{\eps_0}$ and $I_+^{\eps_0}$
do not contain the point $(0,0)$.
 Clearly, $\bigcap_{\eps>0}I_{+}^{\eps}=I_{+}$ and
$\bigcap_{\eps>0}I_{-}^{\eps}=I_{-}$.

Our aim is to show that the generalized solution to the problem
(\ref{eq:1})--(\ref{eq:3}), whose existence is shown in
\cite{KmitHorm}, admits an {\it associated distribution} or a
{\it delta-wave}. The latter means that the family $(u^{\eps})_{\eps>0}$
of solutions to the system with regularized initial and boundary data
\begin{eqnarray}
(\partial_t  + \lambda(x,t)\partial_x) u^{\eps} &=& p(x,t)u^{\eps} + f(x,t,u^{\eps})
  \label{eq:1'}\\
u^{\eps}|_{t=0} &=& a_s^{\eps}+a_r \label{eq:2'} \\
u^{\eps}|_{x=0} &=& (c_s^{\eps}+c_r)\int\limits_{0}^{L}(b_s^{\eps}+b_r)u^{\eps}\,dx
   \label{eq:3'}
\end{eqnarray}
has a weak limit. Here
$$
a_s^{\eps}=a_s\ast \vphi_{\eps},\quad b_s^{\eps}=b_s\ast \vphi_{\eps},\quad c_s^{\eps}=c_s\ast \vphi_{\eps},
$$
where mollifiers $\vphi_{\eps}$ are model delta nets, that is,
$$
\vphi_{\eps}(x)=\frac{1}{\eps}\vphi\Bigl(\frac{x}{\eps}\Bigr)
$$
for an arbitrary fixed $\vphi\in\D(\R)$ with $\int\vphi(x)\,dx=1$. Note that
\begin{equation}\label{eq:eps}
a_s^{\eps}=O\Bigl(\frac{1}{\eps}\Bigr),\quad b_s^{\eps}=O\Bigl(\frac{1}{\eps}\Bigr),\quad c_s^{\eps}
=O\Bigl(\frac{1}{\eps}\Bigr)
\end{equation}
 and
\begin{equation}\label{eq:eps_2}
\begin{array}{c}
\displaystyle
\int\limits_0^L|a_s^{\eps}(x)|\,dx\le C,\quad
\int\limits_0^L|b_s^{\eps}(x)|\,dx\le C,\quad
\int\limits_0^{\infty}
|c_s^{\eps}(t)|\,dt\le C,
\end{array}
\end{equation}
where  $C$ does not depend on $\eps$.
 We will consider mollifiers $\vphi$ with
\begin{equation}\label{eq:+1}
\supp\vphi\subset[-1,1].
\end{equation}
This restriction makes no loss of generality, because
if (\ref{eq:+1}) is not true,
then
$\supp\vphi\subset[-d,d]$ for some $d>0$.
Therefore $\supp\vphi_{\eps}\subset[-d\eps,d\eps]$ and
it is enough to replace $I_+^{\eps}$ by
$I_+^{d\eps}$ to keep all
arguments valid, with the result not depending on  $d$.

It follows from (\ref{eq:+1}) that  for all $\eps>0$
\begin{equation}\label{eq:eps_1}
\begin{array}{c}
\displaystyle
\int\limits_{x_i^*-\eps}^{x_i^*+\eps}a_s^{\eps}(x)\,dx=1,\quad 1\le i\le m;\qquad
\int\limits_{x_j-\eps}^{x_j+\eps}b_s^{\eps}(x)\,dx=1,\quad 1\le j\le k;\\
\displaystyle
\int\limits_{t_p-\eps}^{t_p+\eps}c_s^{\eps}(t)\,dt=1,\quad 1\le p\le l.
\end{array}
\end{equation}


Let $T$ be an arbitrary positive real, $\Pi^T=\{(x,t)\in\Pi\,|\,t<T\}$.
We will show that a delta-wave splits up into the sum $w+v$ of the following
kind.

The function $w$
corresponds to the regular part of the problem.
More specifically, for every $T>0$, the restriction of $w$
to $\Pi^T$
is  the
limit of $w^{\eps}$ in $\Con(\Pi^T)$ as $\eps\to 0$, where $w^{\eps}$
for every fixed
$\eps>0$ is a continuous solution to
the nonlinear problem
\begin{eqnarray}
(\partial_t  + \lambda(x,t)\partial_x) w^{\eps} &=& p(x,t)w^{\eps} + f(x,t,w^{\eps})
  \label{eq:4}\\
w^{\eps}|_{t=0} &=& a_r \label{eq:5}\\
w^{\eps}|_{x=0} &=& c_r\int\limits_{0}^{L}[(b_s^{\eps}+b_r)w^{\eps}+b_r
v^{\eps}]\,dx.
   \label{eq:6}
\end{eqnarray}

The function $v$
corresponds to the singular part of the problem and is  the
limit of  $v^{\eps}$ in $\D'(\Pi)$ as $\eps\to 0$, where $v^{\eps}$ for every
fixed $\eps>0$ is a
 continuous solution to
the linear problem
\begin{eqnarray}
(\partial_t  + \lambda(x,t)\partial_x) v^{\eps} &=& p(x,t)v^{\eps}
  \label{eq:7}\\
v^{\eps}|_{t=0} &=& a_s^{\eps} \label{eq:8}\\
v^{\eps}|_{x=0} &=& c_s^{\eps}\int\limits_{0}^{L}[(b_s^{\eps}
+b_r)w^{\eps}
+b_rv^{\eps}]\,dx
+c_r\int\limits_{0}^{L}b_s^{\eps}
v^{\eps}\,dx.
   \label{eq:9}
\end{eqnarray}


\begin{prop}
For every
$\eps\le\eps_0$ there exist a unique $\Con(\overline{\Pi})$-solution $u^{\eps}$
to the problem
(\ref{eq:1'})--(\ref{eq:3'}), a unique
$\Con(\overline{\Pi})$-solution $w^{\eps}$
to the problem
(\ref{eq:4})--(\ref{eq:6}),
 and
a unique $\Con(\overline{\Pi})$-solution $v^{\eps}$
to the problem
(\ref{eq:7})--(\ref{eq:9}).
\end{prop}
\begin{proof}
For every fixed $\eps>0$, (\ref{eq:1})--(\ref{eq:3})
is a special case of the problem
studied in \cite{KmitHorm}. From the proof of \cite[Theorem 3]{KmitHorm} it follows that,
if
Assumptions 4--7 hold and $\eps$ is so small that $(0,0)\not\in I_+^{\eps}$, this problem
has a unique
$\Con(\overline{\Pi})$-solution $u^{\eps}$.

Fix $\eps>0$.
We consider  (\ref{eq:4})--(\ref{eq:6}) and (\ref{eq:7})--(\ref{eq:9})
simultaneously thereby obtaining an initial-boundary value problem for a system of two hyperbolic equations
with respect to $(w^{\eps},v^{\eps})$.
This is another special case of the problem
studied in \cite{KmitHorm}.
Note that we have zero-order compatibility of (\ref{eq:5}), (\ref{eq:6})
 and of (\ref{eq:8}), (\ref{eq:9}), the former by
Assumption 4 and the latter by Assumption 4
and the fact that $(0,0)\not\in I_+^{\eps}$.
From the proof of \cite[Theorem 3]{KmitHorm} it follows that under
Assumptions 4--7 the  problem has a unique
$\Con(\overline{\Pi})$-solution $(w^{\eps},v^{\eps})$.
\end{proof}

We are now prepared to state  the main result of the paper.
\begin{thm}
Let Assumptions 1--8 hold. Let $u^{\eps}$, for every $\eps>0$,
be the continuous solution to the problem
(\ref{eq:1'})--(\ref{eq:3'}).
Then
$$
u^{\eps}\to w+v \quad {\mbox in}\,\,\,\D'(\Pi) \quad {\mbox as}\,\,\,\,\eps\to 0,
$$
where
\begin{itemize}
\item
for every $T>0$,  the restriction of $w$
to $\Pi^T$
is  the
limit of $w^{\eps}$ in $\Con(\Pi^T)$ as $\eps\to 0$
with
 $w^{\eps}$ being  the continuous solution to
the problem (\ref{eq:4})--(\ref{eq:6}),
\item
$v=\lim\limits_{\eps\to 0}v^{\eps}$ in $\D'(\Pi)$ with $v^{\eps}$ being the
 continuous solution to
the problem
 (\ref{eq:7})--(\ref{eq:9}). Furthermore,
the restriction of $v$
to
$\Pi\setminus I_+$ is identically equal to 0.
\end{itemize}
\end{thm}

\begin{cor}
$
\sing\supp(w+v)=\sing\supp v=\supp v\subset I_+.
$
\end{cor}
This means that $v$ actually represents the purely singular part of the initial
problem. The proof of the corollary is straightforward.

The proof of Theorem 1 consists of five lemmas whose proofs are given in
Sections~4--8.

\begin{lemma}\label{lemma:veps}
Let Assumptions 1, 2, 4--6, and 8 hold and $v^{\eps}$ be as
in Theorem 1. Then
$$
v^{\eps}\to 0
\quad{\mbox pointwise\,\, off}\,\, I_+
\quad{\mbox as}\,\,\eps\to 0.
$$
\end{lemma}

\begin{lemma}\label{lemma:lloc}
Let Assumptions 1--8 hold and $u^{\eps},v^{\eps},$ and $w^{\eps}$ be as
in Theorem 1. Then
$$
u^{\eps}-v^{\eps}-w^{\eps}\to 0 \quad{\mbox in}\,\, \L_{loc}^1(\Pi)
\quad{\mbox as}\,\,\eps\to 0.
$$
\end{lemma}

\begin{lemma}\label{lemma:C}
Let Assumptions 1--8 hold and
$w^{\eps}$ be as
in Theorem 1. Then
$$
w^{\eps} \quad {\mbox converges\,\,in}\,\,\Con(\overline{\Pi^T}) \quad {\mbox as}\,\,\eps\to 0\label{eq:26}
$$
for an arbitrary fixed $T>0$.
\end{lemma}

\begin{lemma}\label{lemma:D'}
Let Assumptions 1--8 hold and
$v^{\eps}$  be as
in Theorem 1. Then
$$
v^{\eps} \quad {\mbox converges\,\,in}\,\,\D'(\Pi)
\quad {\mbox as}\,\,\eps\to 0.\label{eq:25}
$$
\end{lemma}

\begin{lemma}\label{lemma:v=0}
Let Assumptions 1--8 hold,
$v^{\eps}$  be as
in Theorem 1, and
$v=\lim\limits_{\eps\to 0}v^{\eps}$ in
$\D'(\Pi)$. Then $v$ restricted
to
$\Pi\setminus I_+$ is identically
equal to 0.
\end{lemma}

Theorem 1 now follows from the embedding of
$
\L_{loc}^1(\Pi)$ into $\D'(\Pi)$.

\section {Representation of the problems
(\protect\ref{eq:4})--(\protect\ref{eq:6}) and
(\protect\ref{eq:7})--(\protect\ref{eq:9}) in an
integral-operator form}

The problem (\ref{eq:4})--(\ref{eq:6})
is equivalent to the integral-operator equation
\begin{equation}\label{eq:27_0}
w^{\eps}(x,t)=(Rw^{\eps})(x,t)+
\int\limits_{\theta(x,t)}^{t}\Bigl[f(\xi,\tau,w^{\eps})+(p
w^{\eps})(\xi,\tau)\Bigr]\Big|_{\xi=\om(\tau;x,t)}\,d\tau
\end{equation}
and to the corresponding linearized integral-operator equation
\begin{equation}\label{eq:27}
\begin{array}{c}\displaystyle
w^{\eps}(x,t)=(Rw^{\eps})(x,t)\\
\displaystyle
+\int\limits_{\theta(x,t)}^{t}\biggl[f(\xi,\tau,0)+
w^{\eps}\Bigl(\int\limits_0^1
(\nabla_uf)(\xi,\tau,\sigma w^{\eps})\,d\sigma
+p(\xi,\tau)\Bigr)\biggr]\bigg|_{\xi=\om(\tau;x,t)}\,d\tau
\end{array}
\end{equation}
with boundary operator
\begin{equation}\label{eq:28}
(Rw^{\eps})(x,t)=
\cases{a_{r}(\om(0;x,t))
&if
$\theta(x,t)=0$,
 \cr w^{\eps}(0,\theta(x,t))
&if
$\theta(x,t)>0$. \cr}
\end{equation}
Here
$$
\theta(x,t)=\min\limits_{(\om(\tau;x,t),\tau)\in\d\Pi}\tau.
$$
The boundary function $w^{\eps}(0,t)$ is given by  (\ref{eq:6}).

The problem (\ref{eq:7})--(\ref{eq:9})
is equivalent to the integral-operator
equation
\begin{equation}\label{eq:29_1}
v^{\eps}(x,t)
=(Rv^{\eps})(x,t)+
\int\limits_{\theta(x,t)}^{t}
(pv^{\eps})(\om(\tau;x,t),\tau)\,d\tau
\end{equation}
with  boundary operator
\begin{equation}\label{eq:30}
(Rv^{\eps})(x,t)=
\cases{a_{s}^{\eps}(\om(0;x,t))
&if
$\theta(x,t)=0$,
 \cr v^{\eps}(0,\theta(x,t))
&if
$\theta(x,t)>0$. \cr}
\end{equation}
The  function $v^{\eps}(0,t)$ is given by the formula (\ref{eq:9}).
The continuous solution of (\ref{eq:29_1}) can be expressed in the form
\begin{equation}\label{eq:30_1}
v^{\eps}(x,t)=S(x,t)(Rv^{\eps})(x,t)
\end{equation}
with   $\Con(\overline{\Pi})$-function
\begin{equation}\label{eq:S}
\begin{array}{c}\displaystyle
S(x,t)=1+\int\limits_{\theta(x,t)}^tp(\om(\tau;x,t),\tau)\,d\tau\\\displaystyle
+
\int\limits_{\theta(x,t)}^tp(\om(\tau;x,t),\tau)\,d\tau
\int\limits_{\theta(x,t)}^{\tau}p(\om(\tau_1;\om(\tau;x,t),\tau),\tau_1)
\,d\tau_1+\dots.
\end{array}
\end{equation}

\section{Proof of Lemma 1}

By (\ref{eq:30_1}) it suffices to show for every
$(x,t)\in\Pi\setminus I_+$ that, if
$\eps$ is small enough, then
$(Rv^{\eps})(x,t)=0$. If $\theta(x,t)=0$, the latter is true
by
the equality $(Rv^{\eps})(x,t)=a_s^{\eps}(\om(0;x,t))$ and the fact that
$(\om(0;x,t),0)\not\in I_+$.
 Consider the case that $\theta(x,t)>0$.
Since $\theta(x,t)\not\in
\Bigl(\overline{I_+^{\eps}}\cap\{(x,t)\,|\,x=0\}\Bigr)$,
the proof will be complete by showing that
\begin{equation}\label{eq:line}
\begin{array}{c}
\supp v^{\eps}(0,t)\subset
\Bigl(\overline{I_+^{\eps}}\cap\{(x,t)\,|\,x=0\}\Bigr),
\end{array}
\end{equation}
where $v^{\eps}(0,t)$ is defined by (\ref{eq:9}). Observe that
(\ref{eq:line}) is true for the first summand in (\ref{eq:9}). Indeed, by
(\ref{eq:+1}), Assumption 2, and the definition of $I_+^{\eps}$,
\begin{equation}\label{eq:r3}
\begin{array}{c}
\displaystyle
\supp \biggl(c_s^{\eps}\int\limits_{0}^{L}[(b_s^{\eps}
+b_r)w^{\eps}
+b_rv^{\eps}]\,dx\biggr)\subset\bigcup\limits_{i=1}^l[t_i-\eps,t_i+\eps]
\subset\supp c_s^{\eps}\\\displaystyle
\subset
\Bigl(\overline{I_+^{\eps}[0]}\cap\{(x,t)\,|\,x=0\}\Bigr).
\end{array}
\end{equation}
To obtain   (\ref{eq:line}), for the second summand in (\ref{eq:9})
we prove the inclusion
\begin{equation}\label{eq:r4}
\supp \biggl(c_r\int\limits_{0}^{L}b_s^{\eps}
v^{\eps}\,dx\biggr)\subset\bigcup\limits_{n\ge 1}[\tilde t_{i_n}-\eps_{i_n}^-
(\eps),\tilde t_{i_n}+\eps_{i_n}^+
(\eps)].
\end{equation}
Recall that $\tilde t_{i_n}$, $n\ge 1$, are intersection points of
$I_+\setminus I_+[0]$ with the axis $x=0$.
Suppose (\ref{eq:r4}) is false. Then there exists
\begin{equation}\label{eq:r1}
t^1\not\in\bigcup\limits_{n\ge 1}[\tilde t_{i_n}-\eps_{i_n}^-
(\eps),\tilde t_{i_n}+\eps_{i_n}^+
(\eps)]
\end{equation}
such that
\begin{equation}\label{eq:r2}
\int\limits_{0}^{L}b_s^{\eps}(x)
v^{\eps}(x,t^1)\,dx=\int\limits_{0}^{L}
b_s^{\eps}(x)
(Rv^{\eps})(x,t^1)S(x,t^1)\,dx\ne 0.
\end{equation}
We fix such $t^1$ and set
\begin{equation}\label{eq:r5}
J^1=\supp b_s^{\eps}(x)\cap\supp
v^{\eps}(x,t^1).
\end{equation}
By (\ref{eq:r2})
\begin{equation}\label{eq:r6}
\mes J^1\ne 0.
\end{equation}

Assume that $\theta(x_0,t^1)=0$ for some $x_0\in J^1$. By (\ref{eq:30})
and (\ref{eq:r5}),
$$
(Rv^{\eps})(x_0,t^1)=a_{s}^{\eps}
(\om(0;x_0,t^1))\ne 0.
$$
This means that $\om(0;x_0,t^1)\in\supp a_s^{\eps}\subset
\overline{I_+^{\eps}}\cap\{(x,t)\,|\,x=0\}$.
We conclude that
$(x_0,t^1)\in
\overline{I_+^{\eps}[0]}$. Furthermore, from (\ref{eq:r5}) we have
$x_0\in[x_i-\eps,x_i+\eps]$ for some $i\le k$. From the definition
of $I_+^{\eps}$ it follows that, if $(x,t)\in \overline{I_+^{\eps}[j]}$
and $x\in[x_i-\eps,x_i+\eps]$ for some $i\le k$, then $(0,t)\in
\overline{I_+^{\eps}[j+1]}$. Hence $(0,t^1)\in \overline{I_+^{\eps}[1]}$.
This  contradicts
(\ref{eq:r1}).

Assume therefore that $\theta(x,t^1)>0$ for all $x\in J^1$.
Then
in (\ref{eq:r2}) we have
$(Rv^{\eps})(x,t^1)=v^{\eps}(0,\theta(x,t^1))$
and therefore
$$
\int\limits_{0}^{L}b_s^{\eps}(x)
v^{\eps}(0,\theta(x,t^1))S(x,t^1)\,dx\ne 0.
$$
By (\ref{eq:r5}) and (\ref{eq:r6}) there exists $t^2\in\theta(J^1,t^1)$
such that $v^{\eps}(0,t^2)\ne 0$. It is clear that $t^2<t^1$.
Assume that $(0,t^2)\in\overline{I_+^{\eps}}$. Let $x_0$ be such that
$\theta(x_0,t^1)=t^2$. By the definition of $I_+^{\eps}$,
$(x_0,t^1)\in\overline{I_+^{\eps}[j]}$ for some $j\ge 0$. Furthermore,
$x_0\in J^1$ and therefore $x_0\in[x_i-\eps,x_i+\eps]$
for some $i\le k$. Hence $(0,t^1)\in\overline{I_+^{\eps}[j+1]}$.
This again contradicts (\ref{eq:r1}).

Assume therefore that
$(0,t^2)\not\in\overline{I_+^{\eps}}$. On the account of
(\ref{eq:9}) and (\ref{eq:r3}), we rewrite the condition
$v^{\eps}(0,t^2)\ne 0$ as
$$
\int\limits_{0}^{L}b_s^{\eps}(x)
v^{\eps}(x,t^2)\,dx\ne 0.
$$
Set
$$
J^2=\supp b_s^{\eps}(x)\cap\supp v^{\eps}(x,t^2).
$$
Note that
$\mes J^2\ne 0$. Similarly to the above, if
$\theta(x_0,t^2)=0$ for some $x_0\in J^2$, then $(0,t^2)\not\in
\overline{I_+^{\eps}[1]}$, a
contradiction with
(\ref{eq:r1}).
We therefore assume that
$\theta(\xi,t^2)>0$ for all
$x\in J^2$
and continue
in this fashion, thereby constructing sequences $t^k\in\theta(I^{k-1},t^{k-1})$
and $J^k=\supp b_s^{\eps}(x)\cap\supp v^{\eps}(x,t^k)$ for $k\ge 2$
such that $v^{\eps}(0,t^k)\ne 0$ and
$(0,t^k)\not\in \overline{I_+^{\eps}}$. By Assumptions 5 and 6, for some
$$
k\le\biggl\lceil\frac{T\max_{(x,t)\in\overline{\Pi^T}}|\la|}{x_1-\eps_0}\biggr\rceil
$$
there exists $x_0\in J^k$ such that $\theta(x_0,t^k)=0$.
This implies  $(0,t^k)\in\overline{I_+^{\eps}[1]}$,
a contradiction with
(\ref{eq:r1}).

Thus (\ref{eq:r4}) is true and the proof of Lemma 1 is complete.

For the furthure reference observe that
\begin{equation}\label{eq:0}
\supp v^{\eps}\subset\overline{I_+^{\eps}}.
\end{equation}
This fact is true by (\ref{eq:30}), (\ref{eq:30_1}),
(\ref{eq:line}), and
the definition of $I_+^{\eps}$.

\section{Proof of Lemma 2}

Choose $\eps_0$ so small that the number of connected components of
$\Pi^T\cap I_+^{\eps_0}$ and
$\Pi^T\cap I_+$ coincide.
\begin{defn}
Given $T>0$, let
$\Pi_0^T=\{(x,t)
\in\Pi^T\,|\,\om(t;0,0)<x\}$
 and
$\Pi_1^T=\Pi^T\setminus\overline{\Pi_0^T}$.

Let $n(T)$ and $\rho(T)$ be the number of
connected components of    $\Pi_1^T\setminus
\overline{I_+^{\eps}}$ (and $\Pi_1^T\setminus I_+$)
and $\Pi_1^T\cap I_+^{\eps}$ (and $\Pi_1^T\cap I_+$), respectively.
We denote these components, respectively, by
$\Pi^{\eps}(1),$ $\dots,$ $\Pi^{\eps}(n(T))$
($\Pi(1),$ $\dots,$ $\Pi(n(T))$) and
$I_{+}^{\eps}(1),$ $\dots,$
$I_{+}^{\eps}(\rho(T))$ ($I_{+}(1),$ $\dots,$
$I_{+}(\rho(T))$).
\end{defn}

Clearly, $\Pi(i)=
\bigcup_{\eps>0}\Pi^{\eps}(i)$ and
$I_+(i)=\bigcap_{\eps>0}I_{+}^{\eps}(i)$.
   Observe that
$\rho (T)$ does not depend on $\eps$ and either $n(T)=\rho(T)$ or $n(T)=\rho(T)+1$.
 In the latter case, if
$n(T)=\rho(T)+1$, we define
$I_{+}^{\eps}(\rho(T)+1)=\emptyset$.

Given $T$,
we choose $\eps_0$ so small that for all $\eps\le\eps_0$
\begin{equation}\label{eq:supp}
\supp b_s^{\eps}\subset(\om(t_i^*
+\eps_i^+(\eps);0,
t_i^*-\eps_i^-(\eps)),L].
\end{equation}

From (\ref{eq:1'})--(\ref{eq:3'}),
(\ref{eq:4})--(\ref{eq:6}), and (\ref{eq:7})--(\ref{eq:9})
it follows that the difference $u^{\eps}-v^{\eps}-w^{\eps}$ satisfies the system
\begin{eqnarray}
(\partial_t  + \lambda(x,t)\partial_x) (u^{\eps}-v^{\eps}-w^{\eps}) = p(x,t)(u^{\eps}-v^{\eps}-w^{\eps})
\nonumber\\
+F(x,t)(u^{\eps}-v^{\eps}-w^{\eps})+f(x,t,u^{\eps})-f(x,t,u^{\eps}-v^{\eps}),
  \label{eq:10}\\
(u^{\eps}-v^{\eps}-w^{\eps})|_{t=0} = 0, \label{eq:11}\\
(u^{\eps}-v^{\eps}-w^{\eps})|_{x=0} = c_s^{\eps}\int\limits_{0}^{L}b_s^{\eps}
(u^{\eps}-w^{\eps})\,dx\nonumber\\
+(c_s^{\eps}+c_r)\int\limits_{0}^{L}
b_r(u^{\eps}-w^{\eps}-v^{\eps})\,dx
+c_r\int\limits_{0}^{L}
b_s^{\eps}(u^{\eps}-w^{\eps}-v^{\eps})\,dx,
   \label{eq:12}
\end{eqnarray}
where
$$
F(x,t)=\int\limits_0^1
\Bigl(\nabla_uf\Bigr)(x,t,\sigma(u^{\eps}-v^{\eps})+(1-\sigma)w^{\eps})\,d\sigma.
$$


\begin{claim}
 $u^{\eps}-v^{\eps}-w^{\eps}\to 0$ in $\L^1(\overline{\Pi_0^T})$ as $\eps\to 0$.
\end{claim}
\begin{proof}
The problem (\ref{eq:10})--(\ref{eq:12})
on $\Pi_0^T$ reduces to the Cauchy problem (\ref{eq:10})--(\ref{eq:11}).
By Assumption 7, $f$ is globally
bounded  and, by Lemma 1,
$f(x,t,u^{\eps})-f(x,t,u^{\eps}-v^{\eps})\to 0$ as $\eps\to 0$ pointwise off $I_+$.
By Lebesgue's dominated theorem,
$f(x,t,u^{\eps})-f(x,t,u^{\eps}-v^{\eps})\to 0$ in $\L^1
(\overline{\Pi_0^T})$
as $\eps\to 0$.
 Applying Lebesgue's dominated theorem
to the functions $u^{\eps}-v^{\eps}-w^{\eps}$ defined by
(\ref{eq:10})--(\ref{eq:11}),
we obtain the claim.
\end{proof}

\begin{claim}
 $u^{\eps}-v^{\eps}-w^{\eps}\to 0$
pointwise for $(x,t)\in\Pi(1)$ as $\eps\to 0$.
\end{claim}
\begin{proof}
Taking into account Lemma 1, it is sufficient to prove that
$u^{\eps}-w^{\eps}\to 0$  pointwise for $(x,t)\in\Pi(1)$
as $\eps\to 0$. Fix an arbitrary
$(x_0,t_0)\in\Pi(1)$. If $\eps$ is sufficiently small,
the point $(x_0,t_0)$ belongs to $\Pi^{\eps}(1)$, where we have the following
integral representation:
\begin{equation}\label{eq:13}
\begin{array}{c}
\displaystyle
(u^{\eps}-w^{\eps})(x,t)=
c_r(\theta(x,t))\biggl[\int\limits_{0}^{L}b_s^{\eps}(\xi)
(u^{\eps}-w^{\eps})(\xi,\tau)\,d\xi\nonumber
\\\displaystyle
+
\int\limits_{0}^{Q(\tau)}b_r(\xi)(u^{\eps}-w^{\eps})(\xi,\tau)\,d\xi
+\int\limits_{Q(\tau)}^Lb_r(\xi)
(u^{\eps}-w^{\eps}-v^{\eps})(\xi,\tau)\,d\xi\biggr]\bigg|_{\tau=\theta(x,t)}
\\
\displaystyle
+\int\limits_{\theta(x,t)}^t(u^{\eps}-w^{\eps})(\xi,\tau)
\Bigl[p(\xi,\tau)+
\int\limits_{0}^1
(\nabla_uf)(\xi,\tau,\sigma u^{\eps}+(1-\sigma)w^{\eps})\,d\sigma
\Bigr]\Big|_{\xi=\om(\tau;x,t)}\,d\tau.\nonumber
\end{array}
\end{equation}
Here
\begin{equation}\label{eq:Q}
Q(t)=
\cases{\om(t;0,0)
&if
$\theta(L,t)>0$,
 \cr L
&if
$\theta(L,t)=0$. \cr}
\end{equation}
Note that
\begin{equation}\label{eq:w}
\om(t;0,\tau)\le (t-\tau)\max\limits_{(x,t)\in\overline{\Pi^T}}\la(x,t).
\end{equation}
Since, for small enough $\eps$,
$$
|(u^{\eps}-w^{\eps})(x_0,t_0)|\le\max\limits_{(x,t)\in\overline{\Pi^{\eps}(1)}}
|(u^{\eps}-w^{\eps})(x,t)|,
$$
it is sufficient to prove that
$$
\max\limits_{(x,t)\in\overline{\Pi^{\eps}(1)}}
|(u^{\eps}-w^{\eps})(x,t)|=O(\eps).
$$

We start from the evaluation of
the third integral in (\ref{eq:13}). On the account of (\ref{eq:0}),
we represent it in the form
\begin{equation}\label{eq:13_1}
\begin{array}{c}
\displaystyle
\int\limits_{Q(t)}^Lb_r(u^{\eps}-w^{\eps}-v^{\eps})\,dx
=\int\limits_{[Q(t),L]\times\{t\}\setminus(I_+^{\eps}\cap
\{(x,t)\,|\,x\in\R\})}b_r(u^{\eps}-w^{\eps})\,dx\\[10mm]\displaystyle
+\int\limits_{[Q(t),L]\times\{t\}\cap(I_+^{\eps}\cap
\{(x,t)\,|\,x\in\R\})}b_r(u^{\eps}-w^{\eps}-v^{\eps})\,dx.
\end{array}
\end{equation}
If $\theta(L,t)=0$, this integral is equal to 0. Consider the case that
$\theta(L,t)>0$.
To estimate the difference $u^{\eps}-w^{\eps}$ on
$[Q(t),L]\times\{t\}\setminus(I_+^{\eps}
\cap \{(x,t)\,|\,x\in\R\})$,
we consider the corresponding problem
\begin{equation}\label{eq:+}
\begin{array}{c}
\displaystyle
(\partial_t  + \lambda(x,t)\partial_x) (u^{\eps}-w^{\eps}) = \Bigl(p(x,t)+
\int\limits_0^1(\nabla_uf)(x,t,\sigma u^{\eps}+(1-\sigma)w^{\eps})\,d\sigma\Bigr)(u^{\eps}-w^{\eps}),\\
\displaystyle
(u^{\eps}-w^{\eps})|_{t=0} = 0.
\end{array}
\end{equation}
By Assumption 7 this problem has only the trivial solution. Therefore the second integral in
(\ref{eq:13_1}) is equal to 0. We now estimate the third integral.
We have the integral equation
$$
(u^{\eps}-v^{\eps}-w^{\eps})(x,t)
$$
$$
=\int_0^t(p+F)(\xi,\tau)(u^{\eps}-w^{\eps}-v^{\eps})(\xi,\tau)|_{\xi=\om(\tau,x,t)}\,d\tau
$$
$$
+\int_0^t[f(\xi,\tau,u^{\eps})-f(\xi,\tau,u^{\eps}-v^{\eps})]|_{\xi=\om(\tau,x,t)}\,d\tau
$$
that corresponds to the Cauchy problem  (\ref{eq:10})--(\ref{eq:11}).
Combining it with
Assumption~7,
we conclude that
\begin{equation}\label{eq:bound}
\max\limits_{(x,t)\in\overline{\Pi_0^T}\cap\overline{I_+^{\eps}}}
|u^{\eps}-v^{\eps}-w^{\eps}|\le C_1
\end{equation}
for a positive constant $C_1$ not depending on
$\eps$.
Therefore
 the absolute value
 of the third integral in (\ref{eq:13_1}) is bounded from above by
$C_1\eps$.  As a consequence,
\begin{equation}\label{eq:14}
\Bigl|\int\limits_{Q(t)}^Lb_r(u^{\eps}-w^{\eps}-v^{\eps})\,dx\Bigr|
\le C_2\eps.
\end{equation}

In the rest of the proof, $C_i$ for $i\ge 1$
are positive constants that do not depend on
$\eps$.
We will distinguish two cases.

{\it Case 1. $\supp b_s^{\eps_0}\subset[\om(t_1^*;0,0),L]$.}
The first integral in (\ref{eq:13}) vanishes, since
\begin{equation}\label{eq:Pi_1}
\int\limits_{0}^{L}b_s^{\eps}
(u^{\eps}-w^{\eps})\,dx=\int\limits_{[Q(t),L]\times\{t\}\setminus(I_+^{\eps}\cap
\{(x,t)\,|\,x\in\R\})}b_s^{\eps}
(u^{\eps}-w^{\eps})\,dx
\end{equation}
and $u^{\eps}-w^{\eps}$ on $[Q(t),L]\times\{t\}\setminus(I_+^{\eps}\cap
\{(x,t)\,|\,x\in\R\})$ satisfies the problem (\ref{eq:+}) which has only
the trivial solution.
Taking into account (\ref{eq:eps}),
(\ref{eq:13})--(\ref{eq:w}), and (\ref{eq:14}),
similarly to \cite[p. 644]{KmitHorm} we obtain the following
 estimate that holds on $\overline{\Pi^{\eps}(1)}\cap\overline{\Pi^{t^0}}$:
$$
|(u^{\eps}-w^{\eps})(x,t)|\le \frac{C_2\eps}{1-q^0t^0},
$$
where
$$
q^0=
\max\limits_{(x,t,y)\in\overline{\Pi^{T}}\times\R}|(\nabla_uf)(x,t,y)|
$$
$$
+\max\limits_{(x,t)\in
\overline{\Pi^{T}}}|p(x,t)|+
\max\limits_{t\in[0,T]
}|c_r(t)|\max\limits_{x\in[0,L]
}|b_r(x)|
\max\limits_{(x,t)\in
\overline{\Pi^{T}}}|\la(x,t)|,
$$
$$
t^0<q^0.
$$
Iterating this estimate at most
  $\lceil T/t^0\rceil$ times,
 each time using the final
estimate for $|u^{\eps}-w^{\eps}|$
from a preceding iteration,
we  obtain the following
bound that holds on $\overline{\Pi^{\eps}(1)}$:
\begin{equation}\label{eq:16}
|(u^{\eps}-w^{\eps})(x,t)|\le C_3\eps.
\end{equation}
This completes the proof in Case 1.

{\it Case 2. $\supp b_s^{\eps_0}\not\subset[\om(t_1^*;0,0),L]$.}
We fix an arbitrary sequence $0=p_0<p_1<p_2<\dots<p_M=t_1^*$
such that $\supp b_s^{\eps_0}\subset
[\om(p_j;0,p_{j-1}),L]$. Since
$\supp b_s^{\eps}\subset\supp b_s^{\eps_0}$
for $\eps\le\eps_0$, we can choose the same sequence for all
$\eps\le\eps_0$.
Given this sequence,
we devide $\overline{\Pi^{\eps}(1)}$  into a finite number of subsets
\begin{equation}\label{eq:Pi(1,j)}
\Pi^{\eps}(1,j)=\Bigl\{(x,t)\in
\overline{\Pi^{\eps}(1)}\,|\,
\tilde\om(x;0,p_{j-1})\le t\le
\tilde\om(x;0,p_j)\Bigr\}.
\end{equation}
 We prove (\ref{eq:16})
with an appropriate choice of
$C_3$
separately for each of $\Pi^{\eps}(1,j)$. Since
$\supp b_s^{\eps}\subset
[\om(p_1;0,0),L]$, the conditions of Case 1 are true for
$\Pi^{\eps}(1,1)$, and therefore the estimate (\ref{eq:16})
is true for this subset.
The analog of (\ref{eq:16}) for $\Pi^{\eps}(1,2)$ can be obtained
in much the same way. We concentrate only on  changes. For the first
integral in (\ref{eq:13}) we use the representation (\ref{eq:Pi_1})
with one more summand in the right hand side
$$
\int\limits_{\om(t;0,p_1)}^{Q(t)}b_s^{\eps}
(u^{\eps}-w^{\eps})\,dx.
$$
The absolute value of this integral is bounded from above by $C_4\eps$ due to
(\ref{eq:16}) on $\Pi^{\eps}(1,1)$.

To derive (\ref{eq:14}) on
$\Pi^{\eps}(1,2)$ with $\om(t;0,p_1)$ in place of $Q(t)$ and with
new $C_2$, we observe
that  in the analog of
(\ref{eq:13_1}) there appears the third summand
$$
\int\limits_{\om(t;0,p_1)}^{Q(t)}b_r
(u^{\eps}-w^{\eps})\,dx
$$
that can be bounded from above
by
using
 (\ref{eq:16}) for $\Pi^{\eps}(1,1)$.
Similar arguments apply to all subsequent $\Pi^{\eps}(1,j)$. Thus
the estimate (\ref{eq:16}) is true  for
the whole
$\Pi^{\eps}(1)$.

The proof of Claim 2 is complete.
\end{proof}

\begin{claim}
 The functions $u^{\eps}-w^{\eps}-v^{\eps}$ are bounded on
$\overline{I_{+}^{\eps}(1)}$,
uniformly in $\eps$.
\end{claim}
\begin{proof}
Two cases are possible.

{\it Case 1. $(0,t_1)\in I_{+}(1)$.}
We have
\begin{equation}\label{eq:17}
\begin{array}{c}
\displaystyle
(u^{\eps}-v^{\eps}-w^{\eps})|_{x=0} =
G^{\eps}(t)+
(c_s^{\eps}+c_r)\int\limits_{0}^{\om(t;0,t_1-\eps)}
b_r(u^{\eps}-w^{\eps}-v^{\eps})\,dx,\\
\displaystyle
t\in[t_1-\eps,t_1+\eps],
\end{array}
\end{equation}
where
\begin{equation}\label{eq:17_1}
\begin{array}{c}
\displaystyle
G^{\eps}(t)=
(c_s^{\eps}+c_r)\int\limits_{\om(t;0,t_1-\eps)}^{Q(t)}(b_s^{\eps}+b_r)(u^{\eps}-w^{\eps})\,dx
+(c_s^{\eps}+c_r)\int\limits_{Q(t)}^Lb_s^{\eps}(u^{\eps}-w^{\eps})\,dx\\
\displaystyle
+(c_s^{\eps}+c_r)\int\limits_{Q(t)}^{L}
b_r(u^{\eps}-w^{\eps}-v^{\eps})\,dx.
\end{array}
\end{equation}
This representation
follows from (\ref{eq:supp}), (\ref{eq:r3}), and (\ref{eq:r4}).
We now show that $G^{\eps}(t)$ is bounded on $t\in[t_1-\eps,t_1+\eps]$.
Since
$[\om(t;0,t_1-\eps),Q(t)]\times\{t\}$
$\subset$ $\overline{\Pi^{\eps}(1)}$ for
$t\in[t_1-\eps,t_1+\eps]$, the estimate (\ref{eq:16}) on
$\overline{\Pi^{\eps}(1)}$
applies for the  difference $u^{\eps}-w^{\eps}$
under the first integral
in (\ref{eq:17_1}).
Using also
 (\ref{eq:eps}) and (\ref{eq:eps_2}), we conclude
that the first summand in (\ref{eq:17_1}) is bounded uniformly
in $\eps$.
Since
 $u^{\eps}-w^{\eps}\equiv 0$ on
$\Pi_0^T\setminus I_+^{\eps}$
(see (\ref{eq:+})),  the  second integral
is equal to 0.
Applying
(\ref{eq:14}) and (\ref{eq:eps}) to the third summand,
we  see that $G^{\eps}(t)$
are bounded on
$[t_1-\eps,t_1+\eps]$,
uniformly in $\eps$.

Observe  that
\begin{equation}\label{eq:18}
\begin{array}{c}
\displaystyle
(u^{\eps}-v^{\eps}-w^{\eps})(x,t)=(u^{\eps}-v^{\eps}-w^{\eps})(0,\theta(x,t))\\
\displaystyle
+\int\limits_{\theta(x,t)}^{t}p(\xi,\tau)(u^{\eps}-v^{\eps}-w^{\eps})
(\xi,\tau)|_{\xi=\om(\tau;x,t)}\,d\tau\\
\displaystyle
+\int\limits_{\theta(x,t)}^{t}\Bigl[f(\xi,\tau,u^{\eps})-f(\xi,\tau,w^{\eps})
\Bigr]
\Big|_{\xi=\om(\tau;x,t)}\,d\tau,  \quad (x,t)\in \overline{I_{+}^{\eps}(1)},
\end{array}
\end{equation}
where the boundary function $(u^{\eps}-w^{\eps}-v^{\eps})(0,t)$ is given
by (\ref{eq:17}). By Gronwall's
argument applied to $|(u^{\eps}-v^{\eps}-w^{\eps})(x,t)|$,
we easily obtain the estimate
\begin{equation}\label{eq:19}
\begin{array}{c}\displaystyle
|(u^{\eps}-v^{\eps}-w^{\eps})(x,t)|\le C_3\biggl[2T\max
\limits_{(x,t,y)\in\overline{\Pi^T}\times\R}
|f(x,t,y)|+
\max\limits_{(x,t)\in \overline{I_{+}^{\eps}(1)}}\Bigl|\Bigl[G^{\eps}(\tau)\\
\displaystyle
+(c_s^{\eps}+c_r)(\tau)\int\limits_{0}^{\om(\tau;0,t_1-\eps)}
b_r(\xi)(u^{\eps}-v^{\eps}-w^{\eps})(\xi,\tau)\,d\xi\Bigr]\Big|_{\tau=\theta(x,t)}\Bigr|\biggr],
  \quad (x,t)\in
\overline{I_{+}^{\eps}(1)}.
\end{array}
\end{equation}
By (\ref{eq:w}) we have
\begin{equation}\label{eq:20}
\om(t;0,t_1^*-\eps_1^-(\eps))\le C_5\eps
\end{equation}
for $t\in[t_1^*-\eps_1^-(\eps),
t_1^*+\eps_1^+(\eps)]$.
Given a mollifier $\vphi(t)$, let
\begin{equation}\label{eq:q(eps)}
q(\eps)=C_3C_5\max\limits
_{x\in[0,\om(t_1^*+\eps_1^+(\eps);
0,t_1^*-\eps_1^-(\eps))]}|b_r(x)|\Bigl(\max\limits_{t\in[0,T]}|\vphi(t)|+\eps
\max\limits_{t\in[0,T]}|c_r(t)|\Bigr).
\end{equation}
By Assumptions 3 and 5,
\begin{equation}\label{eq:q}
\lim\limits_{\eps\to 0}q(\eps)=0.
\end{equation}
We
choose $\eps$ so small that
\begin{equation}\label{eq:21}
q(\eps)<1.
\end{equation}
On the account of  (\ref{eq:19}),
(\ref{eq:q(eps)}), and (\ref{eq:21}), for sufficiently small $\eps$
we obtain
\begin{equation}\label{eq:21'}
\max\limits_{(x,t)\in\overline{I_{+}^{\eps}(1)}}|(u^{\eps}-v^{\eps}-w^{\eps})(x,t)|
\le \frac{C_3}{1-q(\eps)}\Bigl[2T\max\limits_{(x,t,y)\in\overline{\Pi^T}\times\R}|f(x,t,y)|+
\max\limits_{t\in[t_1-\eps,t_1+\eps]}
|G^{\eps}(t)|\Bigr].
\end{equation}

{\it Case 2. $(0,t_1)\notin I_{+}(1)$.}
By (\ref{eq:supp}), (\ref{eq:r3}), and (\ref{eq:r4})
we have the equality (\ref{eq:17}) with $t_1-\eps$ and $t_1+\eps$ replaced by
$t_1^*-\eps_1^-(\eps)$ and
$t_1^*+\eps_1^+(\eps)$, respectively, and
with
\begin{equation}\label{eq:22}
\begin{array}{c}\displaystyle
G^{\eps}(t)=c_r\int\limits_{\om(t;0,t_1^*-\eps_1^-(\eps))}^{Q(t)}
(b_s^{\eps}+b_r)(u^{\eps}-w^{\eps})\,dx
+c_r\int\limits_{Q(t)}^{L}b_s^{\eps}(u^{\eps}-w^{\eps})
\,dx\\              \displaystyle
+c_r\int\limits_{Q(t)}^{L}b_r(u^{\eps}-w^{\eps}-
v^{\eps})\,dx.
\end{array}
\end{equation}
To estimate the absolute value of the first integral in (\ref{eq:22})
we apply
 (\ref{eq:16}) on $\overline{\Pi^{\eps}(1)}$ and
(\ref{eq:eps_2}). The
second summand is equal to 0
(see  (\ref{eq:+})). For the third integral we use
(\ref{eq:14}).
It follows
that $G^{\eps}(t)$ is bounded on
$[t_1^*-\eps_1^-(\eps),t_1^*+\eps_1^+(\eps)]$,
 uniformly in $\eps$.
The rest of the proof  runs as  in Case 1,
the minor changes being in using  (\ref{eq:19}) and (\ref{eq:21'}) with
$t_1^*-\eps_1^-(\eps)$ in place of
$t_1-\eps$ and $t_1^*+\eps_1^+(\eps)$ in place of
$t_1+\eps$.
\end{proof}

\begin{claim}
1. For every $j\ge 1$, $u^{\eps}-w^{\eps}-v^{\eps}\to 0$  pointwise for $(x,t)\in
\Pi(j)$, as $\eps\to 0$.
2. For every $j\ge 1$, the functions $u^{\eps}-w^{\eps}-v^{\eps}$ are bounded on
$\overline{I_{+}^{\eps}(j)}$,
uniformly in $\eps$.
\end{claim}

\begin{proof}
Items 1 and 2 of the claim follow from the bounds
\begin{equation}\label{eq:C_5}
\begin{array}{c}\displaystyle
\max\limits_{(x,t)\in\overline{\Pi^{\eps}(j)}}|(u^{\eps}-w^{\eps})(x,t)|
\le A_{j}\eps
\end{array}
\end{equation}
and
\begin{equation}\label{eq:C_6}
\begin{array}{c}\displaystyle
\max\limits_{(x,t)\in\overline{I_+^{\eps}(j)}}|(u^{\eps}-v^{\eps}-w^{\eps})(x,t)|
\le B_{j},
\end{array}
\end{equation}
respectively, where $A_j$ and $B_j$ are constants depending only on $j$.
We prove (\ref{eq:C_5}) and
(\ref{eq:C_6}) by induction on $j$.
The base case of $j=1$ is given by
Claims 2 and 3.
Assume that (\ref{eq:C_5}) and (\ref{eq:C_6})
are true for all $j<i$, $i\ge 2$, and prove these estimates for $j=i$.

To prove
(\ref{eq:C_5}) for $j=i$,
we follow the proof of Claim 2 with the following changes.
We use the formula (\ref{eq:13}) with
$\om(\tau;0,t_{i-1}^*+\eps_{i-1}^+(\eps))$ in place of
$Q(\tau)$.
To estimate the third integral in the analog of (\ref{eq:13}), we
represent it in the form
$$
\int\limits_{\om(t;0,t_{i-1}^*+\eps_{i-1}^+(\eps))}^{Q(t)}b_r
(u^{\eps}-w^{\eps}-v^{\eps})\,dx+
\int\limits_{Q(t)}^Lb_r
(u^{\eps}-w^{\eps}-v^{\eps})\,dx,
$$
and  apply the induction assumptions and (\ref{eq:14}).
As a consequence, we obtain the estimate
 (\ref{eq:14})
with $Q(t)$ replaced by $\om(t;0,t_{i-1}^*+\eps_{i-1}^+(\eps))$
and with new $C_2$.
Similarly to Claim 2,
we distinguish two cases.

{\it Case 1. $\supp b_s^{\eps_0}\subset[\om(t_i^*;0,t_{i-1}^*-\eps_{i-1}^-
(\eps_0)),L]$.}
On the account of (\ref{eq:+}), we can rewrite the first summand
in the analog of (\ref{eq:13}) in the form
$$
\int\limits_{0}^{L}b_s^{\eps}
(u^{\eps}-w^{\eps})\,dx=\int\limits_{[\om(t_i^*-\eps_i^-
(\eps);0,t_{i-1}^*-\eps_{i-1}^-
(\eps)),Q(t)]\times\{t\}\setminus(I_+^{\eps}\cap
\{(x,t)\,|\,x\in\R\})}b_s^{\eps}
(u^{\eps}-w^{\eps})\,dx.
$$
Applying (\ref{eq:eps}) and
(\ref{eq:C_5}) for $j<i$, we conclude that the absolute value of the
integral is bounded from above by $C_7\eps$. The rest of the proof for this
case runs similarly to the proof of Claim 2 in Case 1.

{\it Case 2. $\supp b_s^{\eps_0}\not\subset[\om(t_i^*;0,t_{i-1}^*-\eps_{i-1}^-
(\eps_0)),L]$.}
We fix an arbitrary sequence $t_{i-1}^*=p_0<p_1<p_2<\dots<p_M=t_i^*$
such that $\supp b_s^{\eps_0}\subset
[\om(p_j;0,p_{j-1}),L]$.  Given this sequence, we devide
$\overline{\Pi^{\eps}(i)}$  into a finite number of subsets
\begin{equation}
\Pi^{\eps}(i,j)=\Bigl\{(x,t)\in
\overline{\Pi^{\eps}(i)}\,|\,
\tilde\om(x;0,p_{j-1})\le t\le
\tilde\om(x;0,p_j)\Bigr\}.
\end{equation}
Observe that the partition of $\overline{\Pi^{\eps}(i)}$ is finite
for every $\eps>0$ and the number of subsets does not depend on
$\eps$. We further apply arguments similar to those used
in the proof of Claim 2
for Case 2.

To prove
(\ref{eq:C_6}) for $j=i$,
we follow the proof of Claim 3 with the following changes.
Similarly to Claim 3,
we distinguish two cases.

{\it Case 1. $(0,t_k)\in I_{+}(i)$ for some
$k\le l$}.
 We
use the formula
(\ref{eq:17})
with $t_k$ in place of $t_1$ and with
\begin{equation}\label{eq:i38}
\begin{array}{c}
\displaystyle
G^{\eps}(t)=
(c_s^{\eps}+c_r)\int\limits_{[\om(t;0,t_i^*-\eps_i^-(\eps)),L]\times\{t\}
\setminus(I_+^{\eps}\cap\{(x,t)\,|\,x\in\R\})}
(b_s^{\eps}+b_r)(u^{\eps}-w^{\eps})\,dx\\[10mm]
\displaystyle
+(c_s^{\eps}+c_r)\int\limits_
{[\om(t;0,t_i^*-\eps_i^-(\eps)),L]\times\{t\}
\cap(I_+^{\eps}\cap\{(x,t)\,|\,x\in\R\})}
b_r(u^{\eps}-w^{\eps}-v^{\eps})\,dx.
\end{array}
\end{equation}
Estimation of the first summand is based on the inclusion
$$
[\om(t;0,t_i^*-\eps_i^-(\eps)),L]\times\{t\}
\setminus(I_+^{\eps}\cap\{(x,t)\,|\,x\in\R\})\subset
\bigcup_{j=1}^{i}
\overline{\Pi^{\eps}(j)}\cup(\overline{\Pi_0^T}\setminus I_+^{\eps})
$$
and (\ref{eq:C_5}), which for $j<i$ is given by
the induction assumptions and
for $j=i$ is just proved.
Estimation of the second summand is based on the inclusion
$$
[\om(t;0,t_i^*-\eps_i^-(\eps)),L]\times\{t\}
\cap
(I_+^{\eps}\cap\{(x,t)\,|\,x\in\R\})\subset
\bigcup_{j=1}^{i-1}
\overline{I_+^{\eps}(j)}\cup(\overline{\Pi_0^T}\cap
\overline{I_+^{\eps}}),
$$
$$
t\in[t_1-\eps,t_1+\eps],
$$
and (\ref{eq:C_6}) for $j<i$.

{\it Case 2.  $(0,t_k)\notin I_{+}(i)$ for all $k\le l$}.
We use the formula
(\ref{eq:17})
with $t_i^*-\eps_i^-(\eps)$ and $t_i^*+\eps_i^+(\eps)$ in place of
$t_1-\eps$  and $t_1-\eps$, respectively,
and with
\begin{equation}
\begin{array}{c}
\displaystyle
G^{\eps}(t)=
c_r\int\limits_{[\om(t;0,t_i^*-\eps_i^-(\eps)),L]\times\{t\}
\setminus(I_+^{\eps}\cap\{(x,t)\,|\,x\in\R\})}
(b_s^{\eps}+b_r)(u^{\eps}-w^{\eps})\,dx\\[10mm]
\displaystyle
+c_r\int\limits_
{[\om(t;0,t_i^*-\eps_i^-(\eps)),L]\times\{t\}
\cap(I_+^{\eps}\cap\{(x,t)\,|\,x\in\R\})}
(b_s^{\eps}+b_r)(u^{\eps}-w^{\eps}-v^{\eps})\,dx.
\end{array}
\end{equation}

In order to prove the boundedness of $G^{\eps}(t)$ we apply
(\ref{eq:C_5})
for $j\le i$ and (\ref{eq:C_6})   for
$j<i$.

The rest of the proof for both  cases runs similarly to   the
proof of Claim~4 in Case~1.
\end{proof}

From Claim 4, (\ref{eq:C_5}) for
$j\le n(T)$, and (\ref{eq:0})
we conclude that the family
$(u^{\eps}-w^{\eps}-v^{\eps})_{\eps>0}$ is bounded on $\overline{\Pi_1^T}$
uniformly in $\eps$ and converges to 0
almost everywhere in $\overline{\Pi_1^T}$.
By dominated convergence theorem this
family converges to 0
in $\L^1(\Pi_1^T)$-norm.
On the account of Claim 1, $(u^{\eps}-w^{\eps}-v^{\eps})_{\eps>0}$
converges to 0 in $\L^1(\Pi^T)$-norm.
Since
$T$ is arbitrary, this is precisely the assertion of Lemma 2.

\section{Proof of Lemma 3}

Given $T>0$, we choose $\eps_0$ so small that, for all $\eps\le\eps_0$,
the conditions
(\ref{eq:supp}) and
\begin{equation}\label{eq:100}
q(\eps)\exp\Bigl\{T\max\limits_{(x,t)\in\overline{\Pi^T}}|p(x,t)|\Bigr\}<1
\end{equation}
are fulfilled.
Here $q(\eps)$ is defined by
(\ref{eq:q(eps)}). The condition (\ref{eq:100}) follows from
(\ref{eq:q}).

\begin{cla}
The family of functions
$w^{\eps}$ converges
in $\Con(\overline{\Pi_0^T})$ as $\eps\to 0$.
\end{cla}
\begin{proof}
For $w^{\eps}$ on $\overline{\Pi_0^T}$ we use the representation
given by
(\ref{eq:27}) and (\ref{eq:28}).
 Since $(Rw^{\eps})(x,t)$
$=a_r(\om(0;x,t))$
on $\overline{\Pi_0^T}$, the function $w^{\eps}$
for each $\eps>0$ satisfies the same Volterra integral
 equation of the second kind.
This means that $w^{\eps}$
does not depend on $\eps$ and for each
$\eps>0$ is equal to the same continuous function $w(x,t)$ that can be found
from the integral equation
(\ref{eq:27})
by the method of sequential approximation. The claim follows.
\end{proof}

Therewith we are done in $\overline{\Pi_0^T}$. Since
$\overline{\Pi^T}=\overline{\Pi_0^T}\cup\overline{\Pi_1^T}$
and $w^{\eps}$ for each $\eps>0$ is continuous on
$\overline{\Pi^T}$ (see Proposition 1 in Section 2), it
remains to prove the convergence of $w^{\eps}$ in
$\Con(\overline{\Pi_1^T})$. We will check the Cauchy criterion of the
uniform convergence of $w^{\eps}$.
Given $\de>0$, we have to show
for some $\eps_2=\eps_2(\de)$ and every $\eps_1<\eps_2$ that
\begin{equation}\label{eq:*}
\Bigl|\Bigl(w^{\eps_1}-w^{\eps_2}\Bigr)
(x,t)\Bigr|\le\de
\end{equation}
 for all $(x,t)\in\overline{\Pi_1^T}$.

Because of so strong interaction of the regular and the singular parts
(see the problems (\ref{eq:4})--(\ref{eq:6}) and (\ref{eq:7})--(\ref{eq:9})),
in the course of the proof of (\ref{eq:*}) we will need in parallel
to prove  some properties of $v^{\eps}$.

Let
\begin{equation}
D_k(\eps)=
\int\limits_{t_k^*-\eps_k^-(\eps)}
^{t_k^*+\eps_k^+(\eps)}|v^{\eps}(0,t)|\,dt
\end{equation}
and
\begin{equation}
R_k(\eps_1,\eps_2)=
\int\limits_{t_k^*-\eps_k^-(\eps_2)}
^{t_k^*+\eps_k^+(\eps_2)}(v^{\eps_1}-v^{\eps_2})(0,t)\,dt.
\end{equation}

We will prove by induction
on $j$
the following 5 assertions for $1\le j\le n(T)$ ($n(T)$ as well  as
$\Pi(k)$ below are defined by
Definition 3).
Recall that $n(T)$ does not depend on $\eps_2$.
Throughout this section $C$ is a large enough constant that does not
depend on $\eps$.
\vskip1.0mm
{\bf Assertion 1.}
For every $\de>0$, if $\eps_2$ is small enough and $\eps_1<\eps_2$, then
(\ref{eq:*}) is true for all
$\overline{\Pi^{\eps_2}(j)}\cup\overline{I_+^{\eps_2}(j)}$.

\vskip1.0mm
{\bf Assertion 2.}
The functions $w^{\eps}$ are
bounded on
$\overline{\Pi^{\eps}(j)}\cup\overline{I_+^{\eps}(j)}$, uniformly in
$\eps>0$.

\vskip1.0mm
{\bf Assertion 3.}
The estimate $D_j(\eps)\le C$
is true for all $\eps>0$.

\vskip1.0mm
{\bf Assertion 4.}
If $\eps_2$ is small enough and $\eps_1<\eps_2$, then $|R_j(\eps_1,\eps_2)|\le C\eps_2$.

\vskip1.0mm
{\bf Assertion 5.}
$w^{\eps}(x,t)$ converges in $\Con(
\bigcup_{k=1}^j
\overline{\Pi(k)}\cup\overline{\Pi_0^T}\Bigr)$ as $\eps\to 0$.
\vskip1.0mm

Assertion 1 implies the Cauchy criterion of the uniform convergence
of $w^{\eps}$ on $\overline{\Pi_1^T}$. Indeed, given $\de>0$, let
$\eps_2$ be so small that (\ref{eq:*}) is true for every $\eps_1<\eps_2$
on each $\overline{\Pi^{\eps_2}(j)}\cup\overline{I_+^{\eps_2}(j)}$
for $j\le n(T)$. Recall that, for any $\eps_2>0$,
$$
\bigcup\limits_{j=1}^{n(T)}\Bigl(
\overline{\Pi^{\eps_2}(j)}\cup\overline{I_+^{\eps_2}(j)}\Bigr)=
\overline{\Pi_1^T}.
$$
It follows that (\ref{eq:*}) is true on $\overline{\Pi_1^T}$
for all $\eps_1<\eps_2$. By the Cauchy criterion,
$w^{\eps}$ uniformly converges on $\overline{\Pi_1^T}$.

  The
proof of Assertions 1--5 for $j=1$  will be given by Claims 2--10.
The induction step will be carried out
by Claims 11--19.

To prove Assertion 1, we
split $\overline{\Pi^{\eps_2}(j)}\cup\overline{I_+^{\eps_2}(j)}$
into four subsets:
$$
\overline{\Pi^{\eps_2}(j)}\cup\overline{I_+^{\eps_2}(j)}=
\overline{\Pi^{\eps_2}(j)}
\cup
\Bigl(\overline{\Pi^{\eps_1}(j)}\cap\overline{I_+^{\eps_2}(j)}\Bigr)\cup
\overline{I_+^{\eps_1}(j)}\cup
\Bigl(\overline{I_+^{\eps_2}(j)}\cap\overline{\Pi^{\eps_1}(j+1)}\Bigr),
$$
where each two neighboring subsets
have common border. We will prove
Assertion~1 separately for each of the four subsets.


\begin{cla}
The functions $w^{\eps}(x,t)$ are bounded on
$\overline{\Pi^{\eps}(1)}$, uniformly in $\eps$.
\end{cla}
\begin{proof}
We use the representation of $w^{\eps}$ by
(\ref{eq:27_0}) and (\ref{eq:28}) restricted to
$\overline{\Pi^{\eps}(1)}$. In this representation, on the account of
(\ref{eq:6}), (\ref{eq:30}), and
(\ref{eq:30_1}), we have
$$
\Bigl(Rw^{\eps}\Bigr)(x,t)=c_r(\theta(x,t))\int\limits_0^{\om(\theta(x,t);0,0)}
\Bigl(b_s^{\eps}+b_r\Bigr)(\xi)w^{\eps}(\xi,\theta(x,t))\,d\xi
$$
$$
+c_r(\theta(x,t))\int\limits_{\om(\theta(x,t);0,0)}^L
\Bigl[\Bigl(b_s^{\eps}+b_r\Bigr)(\xi)w(\xi,\tau)
+S(\xi,\tau)a_s^{\eps}(\om(0;\xi,\tau))
\Bigr]\Big|_{\tau=\theta(x,t)}\,d\xi,
$$
where $w(x,t)=w^{\eps}(x,t)$ for all
$(x,t)\in\overline{\Pi_0^T}$ and for all $\eps>0$   (see the proof of Claim~1).
Taking into account
(\ref{eq:w}) and the fact that $\theta(x,t)\le t$,
similarly to~\cite[p.~646]{KmitHorm}
we obtain the global  estimate
\begin{equation}\label{eq:103}
\begin{array}{c}\displaystyle
\max\limits_{(x,t)\in\overline{\Pi^{\eps}(1)}}|w^{\eps}(x,t)|
\le\biggl(\frac{1}{1-q^1t^1}\biggr)^
{\lceil\frac{T}{t^1}\rceil} P(E)
\biggl(1+\max\limits_{x\in[0,L]}|a_r(x)|\\
\displaystyle
+\max\limits_{(x,t)\in\overline{\Pi^T}}
|S(x,t)|
\max\limits_{\eps}\int\limits_0^L|a_s^{\eps}(x)|\,dx\biggr),
\end{array}
\end{equation}
where
\begin{equation}\label{eq:110}
\begin{array}{c}\displaystyle
E=\max\limits_{t\in[0,T]}|c_r(t)|\biggl(\max\limits_{x\in[0,L]}|b_r(x)|
+\max\limits_{\eps}\int\limits_0^L|b_s^{\eps}(x)|\,dx\biggr),
\\\displaystyle
q^1=
(1+LE)\biggl(\max\limits_{(x,t,y)\in\overline{\Pi^{T}}\times\R}|(\nabla_uf)(x,t,y)|
+\max\limits_{(x,t)\in
\overline{\Pi^{T}}}|p(x,t)|\biggr)+E\max\limits_{(x,t)\in
\overline{\Pi^{T}}}|\la(x,t)|,
\end{array}
\end{equation}
$t^1$ is  a real so small that
\begin{equation}\label{eq:t_0}
t^1<q^1
\end{equation}
and
\begin{equation}\label{eq:48_1}
\supp b_s^{\eps}\subset(\om(kt^1;0,
(k-1)t^1),L]\quad
\mbox{for\,\, all}\quad 1\le k\le \lceil T/t^1\rceil,
\end{equation}
and $P(E)$ is a polynomial  of degree $\lceil T/t^1\rceil$ with positive
coefficients depending on $f(x,t,0),L$, and $T$.
The claim now follows by
(\ref{eq:eps_2}) and
Assumptions 5 and 7.
\end{proof}

\begin{cla}
1. Provided $\eps_2$ is small enough,
for all $\eps'_2\le\eps_2$ and
for all $\eps'_1\le\eps'_2$ the estimate
\begin{equation}\label{eq:r7}
\Bigl|\Bigl(w^{\eps'_1}-w^{\eps'_2}\Bigr)
(x,t)\Bigr|\le
\de
\end{equation}
is true on
$\overline{\Pi^{\eps_2}(1)}$.\\
2. Provided $\eps_2$ is small enough,
(\ref{eq:*}) is true on $\overline{\Pi^{\eps_2}(1)}$.
\end{cla}
\begin{proof}
 Recall that $\overline{\Pi^{\eps_2}(1)}\subset$
$\overline{\Pi^{\eps'_1}(1)}$ and $\overline{\Pi^{\eps_2}(1)}\subset$
$\overline{\Pi^{\eps'_2}(1)}$. To
represent
 $w^{\eps'_1}$ and $w^{\eps'_2}$ on $\overline{\Pi^{\eps_2}(1)}$, we
will use the system
(\ref{eq:4})--(\ref{eq:6}) restricted to $\overline{\Pi^{\eps_2}(1)}$.
For the difference $w^{\eps'_1}-w^{\eps'_2}$ we will employ
 the corresponding linearized
integral-operator equation.
We distinguish two
cases.

{\it Case 1.
$\supp b_s^{\eps_2}\subset[\om(t_1^*;0,0),L]$.}
Using the fact that
$w^{\eps}(x,t)\equiv w(x,t)$ on $\overline{\Pi_0^T}$ for $\eps>0$,
we obtain the integral equation
\begin{equation}\label{eq:200}
\begin{array}{c}\displaystyle
\Bigl(w^{\eps'_1}-w^{\eps'_2}\Bigr)(x,t)=c_r(\theta(x,t))
\int\limits_0^{\om(\theta(x,t);0,0)}b_r(\xi)(w^{\eps'_1}-w^{\eps'_2})
(\xi,\theta(x,t))\,d\xi\\
\displaystyle
+S_1(x,t)+S_2(\theta(x,t))+
S_3(\theta(x,t)),
\end{array}
\end{equation}
where
\begin{equation}\label{eq:51_i}
\begin{array}{c}\displaystyle
S_1(x,t)=\int\limits_{\theta(x,t)}^t
\biggl[\Bigl[p(\xi,\tau)\\
\displaystyle
+
\int\limits_0^1(\nabla_u f)(\xi,\tau,\sigma w^{\eps'_1}
+(1-\sigma)w^{\eps'_2})\,d\sigma\Bigr]
(w^{\eps'_1}-w^{\eps'_2})(\xi,\tau)\biggr]
\bigg|_{\xi=\om(\tau;x,t)}\,d\tau,\\
\displaystyle
S_2(t)=c_r(t)
\int\limits_{Q(t)}^L(b_s^{\eps'_1}-b_s^{\eps'_2})(x)
w(x,t)\,dx,\\
\displaystyle
S_3(t)=c_r(t)
\int\limits_{Q(t)}^Lb_r(x)
\Bigl(v^{\eps'_1}-v^{\eps'_2}\Bigr)(x,t)\,dx,
\end{array}
\end{equation}
and $Q(t)$ is defined by (\ref{eq:Q}).
We now estimate $|S_2(t)|$ and $|S_3(t)|$.
Since the function $w$ on $\overline{\Pi_0^T}$ is uniformly continuous,
 the properties (\ref{eq:eps_1}) and (\ref{eq:eps_2})
hold, and
$\supp b_s^{\eps}\subset\bigcup_{j=1}^k[x_j-\eps,x_j+\eps]$,
we have
\begin{equation}\label{eq:201}
\begin{array}{c}\displaystyle
|S_2(t)|\le\max\limits_{t\in[0,T]}|c_r(t)|
\int\limits_{Q(t)}^L\Bigl(|b_s^{\eps'_1}|+|b_s^{\eps'_2}|\Bigr)
\sum\limits_{j=1}^k\Bigl|w(x,t)-
\chi_{[x_j-\eps'_2,x_j+\eps'_2]}(x)
w(x_j,t)\Bigr|\,dx\\
\displaystyle
+\max\limits_{t\in[0,T]}|c_r(t)|\sum\limits_{j=1}^k
\max\limits_{(x_j,t)\in\overline{\Pi_0^T}}
|w(x_j,t)|\biggl|\int
\limits_{x_j-\eps'_2}^{x_j+\eps'_2}\Bigl(b_s^{\eps'_1}-b_s^{\eps'_2}
\Bigr)(x)
\,dx\biggr|\le C\eps_2.
\end{array}
\end{equation}
Here $\chi_{\Omega}(x,t)$ denotes the characteristic function
 of a set $\Omega$.

Taking into account (\ref{eq:30_1})
and changing coordinates
$(x,t)$
to $(\om(0;x,t),t)$, we estimate
$|S_3(t)|$ in the following way:
\begin{equation}\label{eq:201_1}
\begin{array}{c}\displaystyle
|S_3(t)|=\biggl|c_r(t)
\int\limits_{Q(t)}^L
\Bigl(b_rS\Bigr)(x,t)
\Bigl(a_s^{\eps'_1}-a_s^{\eps'_2}\Bigr)(\om(0;x,t))\,dx\biggr|\\
\displaystyle
=\biggl|c_r(t)
\int\limits_0^{\om(0;L,t)}
\frac{(b_rS)(\xi,t)}
{(\d_{x}\om)(0;\xi,t)}\Big|_{\xi=\om(t;x,0)}\Bigl(a_s^{\eps'_1}
-a_s^{\eps'_2}\Bigr)
(x)\,dx\biggr|\\
\displaystyle
\le\max\limits_{t\in[0,T]}|c_r(t)|
\int\limits_0^{x_{m(t)}^*+\eps'_2}\Bigl(
|a_s^{\eps'_1}|+|a_s^{\eps'_2}|\Bigr)\\
\displaystyle
\times\sum\limits_{j=1}^{m(t)}\biggl|\frac{(b_rS)(\om(t;x,0),t)}
{(\d_{x}\om)(0;\om(t;x,0),t)}-
\chi_{[x_j^*-\eps'_2,x_j^*+\eps'_2]}(x)
\frac{(b_rS)(\om(t;x_j^*,0),t)}
{(\d_{x}\om)(0;\om(t;x_j^*,0),t)}\biggr|\,dx
\\
\displaystyle
+\max\limits_{t\in[0,T]}|c_r(t)|
\sum\limits_{j=1}^{m(t)}
\max\limits_{(\om(t;x_j^*,0),t)\in\overline{\Pi_0^T}}
\biggl|\frac{(b_rS)(\om(t;x_j^*,0),t)}
{(\d_x\om)(0,\om(t;x_j^*,0),t)}\biggr|
\biggl|
\int\limits_0^{x_{m(t)}^*+\eps'_2}\Bigl(a_s^{\eps'_1}-a_s^{\eps'_2}
\Bigr)(x)\,dx
\biggr|
\\
\displaystyle
+
\max\limits_{t\in[0,T]}|c_r(t)|
\int\limits_{x_{m(t)}^*+\eps'_2}^{\om(0;L,t)}\Bigl(
|a_s^{\eps'_1}|+|a_s^{\eps'_2}|\Bigr)
\biggl|\frac{(b_rS)(\om(t;x,0),t)}
{(\d_{x}\om)(0;\om(t;x,0),t)}
\biggr|\,dx,\end{array}
\end{equation}
where $m(t)$ is the number of indices $j\le m$ such that
$x_j^*+\eps'_2\in[0,\om(0;L,t)].$
Similarly to ~\cite[p.~644]{KmitHorm} we obtain the
estimate
  for $|(w^{\eps'_1}-w^{\eps'_2})(x,t)|$
on $\overline{\Pi^{\eps_2}(1)}\cap\overline{\Pi^{t^1}}$ :
\begin{equation}\label{eq:400}
\Bigl|\Bigl(w^{\eps'_1}-w^{\eps'_2}
\Bigr)(x,t)\Bigr|\le
\frac{C\eps_2}{1-q^1t^1},
\end{equation}
where $q^1$ and $t^1$ are defined
by (\ref{eq:110}) and (\ref{eq:t_0}).
 Indeed, the second summand in the right-hand side of (\ref{eq:201_1})
is equal to 0
by (\ref{eq:eps_1}). To estimate the first summand,
 we use
(\ref{eq:eps_2}) and
the uniform continuity property for
$b_r$, $S$, and $\la$ on $\overline\Pi^T$. To
estimate
the third summand,
we observe that
the integral is equal to 0 if
$\om(0;L,t)\le x_{m(t)}^*+\eps'_2$
and is actually from $x_{m(t)+1}^*-
\eps'_2$ to $\om(0;L,t)$.
In the latter case
$
\om(0;L,t)-x_{m(t)+1}^*+
\eps'_2\le C\eps'_2.
$
Combining this bound with
the continuity of $\la$ and the condition $\om(t;\om(0;L,t),0)=L$,
we obtain
$L-\om(t;x_{m(t)+1}^*-\eps'_2,0)\le C\eps'_2.$
 Since $b_r(L)=0$ by Assumption 3, we conclude that
$$
\max\limits_{x\in[x_{m(t)+1}^*-\eps'_2,\om(0;L,t)]}\Bigl|\Bigl
(b_rS\Bigr)(\om(t;x,0),t)\Bigr|=
\max\limits_{x
\in[\om(t;x_{m(t)+1}^*-\eps'_2,0),L]}\Bigl|\Bigl
(b_rS\Bigr)(x,t)\Bigr|
$$
$$
\le
\max\limits_{(x,t)
\in[L-C\eps'_2,L]\times[0,T]}\Bigl|\Bigl
(b_rS\Bigr)(x,t)\Bigr|
\le C\eps_2.
$$
It follows that
\begin{equation}\label{eq:202}
|S_3(t)|\le C\eps_2.
\end{equation}

Using (\ref{eq:200}), (\ref{eq:201}), and (\ref{eq:202}), we
derive (\ref{eq:400}) by Gronwall's argument applied to
$|w^{\eps'_1}
-w^{\eps'_2}|$.
Iterating this  estimate at most
$\lceil T/t^1\rceil$ times,
each time using the final
estimate for $|(w^{\eps'_1}
-w^{\eps'_2})(x,t)|$ from a preceding
iteration,
we obtain the bound
\begin{equation}\label{eq:203}
\Bigl|\Bigl(w^{\eps'_1}-w^{\eps'_2}\Bigr)(x,t)\Bigr|
\le\biggl(\frac{1}
{1-q^1t^1}\biggr)
^{\lceil\frac{T}{t^1}\rceil}P(E_1)C\eps_2,
\end{equation}
where
$$
E_1=\max\limits_{(x,t)\in\overline{\Pi^{T}}}|c_r(t)b_r(x)|
$$
and $P(E_1)$ is a polynomial of degree
$\lceil T/t^1\rceil$ with positive coefficients depending on
$L$ and $T$.

{\it Case 2.
$\supp b_s^{\eps_2}\not\subset[\om(t_1^*;0,0),L]$.}
We devide $\overline{\Pi^{\eps_2}(1)}$  into a finite number of subsets
$
\Pi^{\eps_2}(1,j)$, $j\le M$, defined by  (\ref{eq:Pi(1,j)})
with $\eps$ replaced by $\eps_2$.
Note that, if $j<M$, then $\Pi^{\eps_2}(1,j)$ actually does not depend
on $\eps_2$.
 We prove an analog of (\ref{eq:203}) with an appropriate choice of
$t^1$, $P$, and $C$,
separately for each of
$\Pi^{\eps_2}(1,j)$. Since
$\supp b_s^{\eps_2}\subset
[\om(p_1;0,0),L]$, the conditions of Case 1 are true for
$\Pi^{\eps_2}(1,1)$, and therefore the estimate (\ref{eq:203})
is true for this subset. Thus, provided $\eps_2$ is small enough,
the estimate (\ref{eq:r7})
holds on $\Pi^{\eps_2}(1,1)$.

The analog of  (\ref{eq:203}) for  $\Pi^{\eps_2}(1,2)$
can be obtained in much the same way. We concentrate only on changes.
 On $\Pi^{\eps_2}(1,2)$ we use the integral equation
\begin{equation}
\begin{array}{c}\displaystyle
\Bigl(w^{\eps'_1}-w^{\eps'_2}\Bigr)(x,t)=c_r(\theta(x,t))
\int\limits_0^{\om(\theta(x,t);0,p_1)}b_r(\xi)(w^{\eps'_1}-w^{\eps'_2})
(\xi,\theta(x,t))\,d\xi\\
\displaystyle
+S_1(x,t)+S_2(\theta(x,t))+
S_3(\theta(x,t))+S_4(\theta(x,t)),
\end{array}
\end{equation}
where
$$
S_4(t)=c_r(t)\int\limits_{\om(t;0,p_1)}^{Q(t)}
\Bigl[b_s^{\eps'_1}(x)
\Bigl(w^{\eps'_1}-
w^{\eps'_2}\Bigr)(x,t)+
\Bigl(b_s^{\eps'_1}-b_s^{\eps'_2}\Bigr)(x)w^{\eps'_2}(x,t)\Bigr]\,dx.
$$
We now bound  $|S_2(t)+S_4(t)|$.
Observe that
$[\om(t;0,p_1),Q(t)]\times\{t\}$
$\subset$ $\Pi^{\eps_2}(1,1)$
if $t\in[p_1,p_2]$.
By Proposition 1, Claim 1,  and the estimate (\ref{eq:r7}) on
$\Pi^{\eps_2}(1,1)$, we conclude that $w^{\eps}$ converges
in $\Con(\Pi^{\eps_2}(1,1)\cup
\overline{\Pi_0^T})$ to a continuous function $w(x,t)$.
Using the equality
$w^{\eps'_1}(x,t)=w^{\eps'_2}(x,t)=w(x,t)$ on $\overline{\Pi_0^T}$,
similarly to (\ref{eq:201}) we derive the bound
$$
|S_2(t)+S_4(t)|
\le\max_{t\in[0,T]}|c_r(t)|\max_{(x,t)\in\Pi(1,1)}
|w^{\eps'_1}-w^{\eps'_2}|\int\limits_{\om(t;0,p_1)}^{Q(t)}
|b_s^{\eps'_1}|\,dx
$$
$$
+\max_{t\in[0,T]}|c_r(t)|\int\limits_{\om(t;0,p_1)}^{L}
\Bigl(|b_s^{\eps'_1}|+|b_s^{\eps'_2}|
\Bigr)\biggl|
w^{\eps'_2}(x,t)-\sum\limits_{j=1}^{k}
\chi_{[x_j-\eps'_2,
x_j+\eps'_2]}(x)w(x_j,t)\biggr|\,dx
$$
$$
+\max_{t\in[0,T]}|c_r(t)|\sum
\limits_{j=1}^{k}
\max\limits_{(x_j,t)\in\overline{\Pi_0^T}\cup\Pi(1,1)}|w(x_j,t)|
\biggl|\int\limits_{x_j-\eps'_2}^{x_j+\eps'_2}\Bigl(b_s^{\eps'_1}-
b_s^{\eps'_2}\Bigr)(x)\,dx
\biggr|
\le C\eps_2.
$$
Now, using (\ref{eq:202}),
we conclude that   (\ref{eq:203})
holds for
$\Pi^{\eps_2}(1,2)$ with new
$t^1$, $P$, and $C$.

Similar arguments apply to
the subsets $\Pi^{\eps_2}(1,j)$. Thus the estimate
\begin{equation}\label{eq:204}
\Bigl|\Bigl(w^{\eps'_1}-w^{\eps'_2}\Bigr)(x,t)\Bigr|\le C\eps_2
\end{equation}
is true for the whole
$\overline{\Pi^{\eps_2}(1)}$ in both cases.

The estimate (\ref{eq:r7}) follows
from  (\ref{eq:204}),
where $\eps_2$ is chosen small enough.
The proof of Item 1 is complete.

Item 2 is a straightforward consequence of Item 1.
\end{proof}

\begin{cla}
Provided $\eps_2$ is small enough,
(\ref{eq:*}) is true on $\overline{\Pi^{\eps_1}(1)}
\cap \overline{I_+^{\eps_2}(1)}$.
\end{cla}
\begin{proof}
The functions $w^{\eps_1}(x,t)$ and $w^{\eps_2}(x,t)$
on $\overline{\Pi^{\eps_1}(1)}
\cap \overline{I_+^{\eps_2}(1)}$ are represented
by (\ref{eq:27_0}) restricted, respectively,
 to $\overline{\Pi^{\eps_1}(1)}$ and
$\overline{I_+^{\eps_2}(1)}$.

Proposition 1 together with Claims 1  and 3 (Item 1)
imply, for each fixed $\eps_2>0$ as small as in Claim~3,
the Cauchy criterion of the uniform convergece of $w^{\eps}$ on
$\overline{\Pi_0^T}\cup
\overline{\Pi^{\eps_2}(1)}$.

 Therefore
\begin{equation}\label{eq:55_1}
w^{\eps}(x,t) \,\,\mbox{converges\,\,in} \,\,
\Con\Bigl(\overline{\Pi_0^T}\cup
\overline{\Pi^{\eps_2}(1)}\Bigr) \,\,\mbox{as}\,\, \eps\to 0.
\end{equation}
As above, the
limit function will be denoted by $w(x,t)$.
Now, using
(\ref{eq:supp}),
we have the  representation
$$
\Bigl(w^{\eps_1}-w^{\eps_2}\Bigr)(x,t)
$$
$$
=\biggl[ c_r(\tau)
\int\limits_{\om(\tau;0,t_1^*-
\eps_1^-(\eps_2))}^L
\Bigl((b_s^{\eps_1}+b_r)
(w^{\eps_1}-w^{\eps_2})\Bigr)(\xi,\tau)\,d\xi\biggr]\bigg|_{\tau=\theta(x,t)}
$$
$$
+\biggl[ c_r(\tau)
\int\limits_{\om(\tau;0,t_1^*-
\eps_1^-(\eps_2))}^L
\Bigl(b_s^{\eps_1}-b_s^{\eps_2}\Bigr)(\xi)
\Bigl(w^{\eps_2}(x,t)
$$
\begin{equation}\label{eq:206}
\begin{array}{c}
\displaystyle
-
\sum\limits_{j=1}^k
\chi_{[x_j-\eps_2,x_j+\eps_2]}(x)w(x_j,t)\Bigr)\,d\xi\biggr]\bigg|_{\tau=\theta(x,t)}
\\
\displaystyle
+\biggl[ c_r(\tau)
\sum\limits_{j=1}^k w(x_j,t)\int\limits_{\om(\tau;0,t_1^*-
\eps_1^-(\eps_2))}^L \Bigl(b_s^{\eps_1}-b_s^{\eps_2}\Bigr)(\xi)\,d\xi
\biggr]\bigg|_{\tau=\theta(x,t)}+S_3(\theta(x,t))
\\
\displaystyle
+\biggl[ c_r(\tau)
\int\limits_0^{\om(\tau;0,t_1^*-
\eps_1^-(\eps_2))}
\Bigl(b_r
(w^{\eps_1}-w^{\eps_2})\Bigr)(\xi,\tau)\,d\xi\biggr]\bigg|_{\tau=\theta(x,t)}
+S_1(x,t)+S_5(\theta(x,t)),
\end{array}
\end{equation}
where
$$
S_5(t)=c_r(t)
\int\limits_0^{\om(t;0,t_1^*-\eps_1^-(\eps_2))}
\Bigl(b_rv^{\eps_2}\Bigr)
(x,t)\,dx.
$$
For the absolute values
of the first four summands in (\ref{eq:206}) we obtain the upper bound
$C\eps_2$ by
the following argument.
Note that  $[\om(t;0,t_1^*-
\eps_1^-(\eps_2)),L]\times\{t\}$ $\subset$
$\overline{\Pi_0^T}\cup
\overline{\Pi^{\eps_2}(1)}$ for
$t\in[t_1^*-\eps_1^-(\eps_2),t_1^*-\eps_1^-(\eps_1)]$.
For the first summand the bound now follows from
(\ref{eq:203}),
(\ref{eq:eps_2}), and Claims 1 and 3 (Item 2).
For the second summand we apply
(\ref{eq:55_1})  and (\ref{eq:eps_2}). For the third one we use the properties
(\ref{eq:eps_1}).
 For $|S_3(t)|$ we use
 the estimate (\ref{eq:202}).

To prove the upper bound $C\eps_2$
for $|S_5(t)|$ we need an estimate for $v^{\eps}$
on $\overline{I_+^{\eps}(1)}$.
To obtain it we consider two cases.

{\it Case 1. $(0,t_1)\in I_{+}(1)$.}
Taking into account (\ref{eq:supp}), we represent $w^{\eps}$ and
$v^{\eps}$
on $\overline{I_{+}^{\eps}(1)}$ in the form
\begin{equation}\label{eq:34_1}
\begin{array}{c}\displaystyle
w^{\eps}(x,t)=\Bigl[c_r(S_6^{\eps}+S_7^{\eps})
\Bigr](\theta(x,t))\\\displaystyle
+\int\limits_{\theta(x,t)}^tf(\om(\tau;x,t),\tau,w^{\eps})\,d\tau
+\int\limits_{\theta(x,t)}^tp(\om(\tau;x,t),\tau)w^{\eps}\,d\tau
\end{array}
\end{equation}
and
\begin{equation}\label{eq:34_2}
\begin{array}{c}
\displaystyle
v^{\eps}(x,t)=\Bigl[c_s(S_6^{\eps}+S_7^{\eps})
\Bigr](\theta(x,t))
+\int\limits_{\theta(x,t)}^tp(\om(\tau;x,t),\tau)v^{\eps}\,d\tau,
\end{array}
\end{equation}
where
\begin{equation}\label{eq:SS}
\begin{array}{c}\displaystyle
S_6^{\eps}(t)=
\int\limits_{\om(t;0,t_1^*-\eps_1^-(\eps))}^L
\Bigl(b_{s}^{\eps}+b_r\Bigr)(x)w^{\eps}(x,t)\,dx+
\int\limits_{Q(t)}^L
\Bigl(b_{r}v^{\eps}\Bigr)(x,t)\,dx,
\\\displaystyle
S_7^{\eps}(t)=
\int\limits_0^{\om(t;0,t_1^*-\eps_1^-(\eps))}b_r(x)(w^{\eps}+v^{\eps})(x,t)\,dx.
\end{array}
\end{equation}
Summing up, we have
\begin{equation}\label{eq:35}
\begin{array}{c}\displaystyle
w^{\eps}(x,t)+v^{\eps}(x,t)=\Bigl[(c_r+c_s^{\eps})
S_6^{\eps}\Bigr](\theta(x,t))+\Bigl[(c_r+c_s^{\eps})
S_7^{\eps}\Bigr](\theta(x,t))\\
\displaystyle
+\int\limits_{\theta(x,t)}^tf(\om(\tau;x,t),\tau,w^{\eps})\,d\tau+
\int\limits_{\theta(x,t)}^tp(\om(\tau;x,t),\tau)(w^{\eps}+v^{\eps})
\,d\tau.
\end{array}
\end{equation}
By (\ref{eq:30}), (\ref{eq:30_1}),  (\ref{eq:eps_2}),
Proposition 1,
  and Claim 2,
 $S_6^{\eps}(t)$
is a continuous function
 and satisfies the uniform in $\eps$ estimate
\begin{equation}\label{eq:S_6}
\begin{array}{c}\displaystyle
|S_6^{\eps}(t)|\le C,\quad t\in[t_1^*-\eps_1^-(\eps),t_1^*+\eps_1^+(\eps)].
\end{array}
\end{equation}
 By Proposition 1, $S_7^{\eps}(t)$
 is   continuous.
We now derive an upper bound for $|S_7^{\eps}(t)|$.
Applying the method of sequential approximation
to the function $w^{\eps}+v^{\eps}$ given by the formula (\ref{eq:35}),
we obtain the estimate
\begin{equation}
\begin{array}{c}\displaystyle
\max\limits_{(x,t)\in
\overline{I_{+}^{\eps}}(1)}|(w^{\eps}+v^{\eps})(x,t)|\le
\biggl[T\max\limits_{(x,t,y)\in\overline{\Pi^{T}}\times\R}|f(x,t,y)|
\nonumber\\
\displaystyle
+\max\limits_{t\in[t_1-\eps,t_1+\eps]}|c_r(t)+c_s^{\eps}(t)|
\Bigl(C+\max\limits_{t\in[t_1-\eps,t_1+\eps]}
\Bigl|S_7^{\eps}(t)|\Bigr)\biggr]
\exp\Bigl\{T\max\limits_{(x,t)\in
\overline{\Pi^T}}|p(x,t)|\Bigr\}.\nonumber
\end{array}
\end{equation}
 By (\ref{eq:100}),
\begin{equation}\label{eq:105}
\begin{array}{c}\displaystyle
\max\limits_{(x,t)\in
\overline{I_{+}^{\eps}(1)}}|w^{\eps}+v^{\eps}|\\
\displaystyle
\le\frac{C}{1-
q(\eps)C\exp\Bigl\{T\max\limits_{(x,t)\in\overline{\Pi^T}}|p|\Bigr\}}
\biggl[\max\limits_{(x,t,y)\in\overline{\Pi^{T}}\times\R}|f(x,t,y)|
+\max\limits_{t\in[t_1-\eps,t_1+\eps]}|c_r(t)+c_s^{\eps}(t)|
\biggr].
\end{array}
\end{equation}
From (\ref{eq:20}), (\ref{eq:105}),
 (\ref{eq:eps}), and Assumption 3 we conclude that
\begin{equation}\label{eq:61_1}
\begin{array}{c}
\displaystyle
|S_7^{\eps}(t)|\le C\max\limits_{x\in[0,\om(t;0,t_1-\eps)]}|b_r(x)|
\le C\max\limits_{x\in[0,C_5\eps]}|b_r(x)|\le C\eps,\\
\displaystyle
t\in[t_1^*-\eps_1^-(\eps),t_1^*+\eps_1^+(\eps)].
\end{array}
\end{equation}
Combining  (\ref{eq:29_1}) and (\ref{eq:30_1}) with
(\ref{eq:34_2}),
we obtain
\begin{equation}\label{eq:205}
v^{\eps}(x,t)=c_s^{\eps}(\theta(x,t))\Bigl(S_6^{\eps}+S_7^{\eps}\Bigr)(\theta(x,t))S(x,t),
\end{equation}
where the functions $|(S_6^{\eps}+S_7^{\eps})(\theta(x,t))|$
 on
$\overline{I_{+}^{\eps}(1)}$
are
bounded uniformly in $\eps$. The latter is true
by (\ref{eq:S_6}) and (\ref{eq:61_1}). The formula
(\ref{eq:205}) implies
\begin{equation}\label{eq:208}
v^{\eps}=O\Bigl(\frac{1}{\eps}\Bigr)
\end{equation}
for $(x,t)\in \overline{I_+^{\eps}(1)}$.
  Taking into account (\ref{eq:205}),
(\ref{eq:20}), (\ref{eq:eps}), and Assumption 3, we  derive the bound
\begin{equation}\label{eq:209}
|S_5(t)|\le C\max\limits_{x\in[0,
\om(t;0,t_1-\eps_2)]}|b_r(x)|
\le C\max\limits_{x\in[0,C_5\eps_2]}|b_r(x)|\le C\eps_2.
\end{equation}

{\it Case 2. $(0,t_1)\not\in I_+(1)$.} The proof is much the same as for Case 1.
The only difference is in evaluation of $|S_5(t)|$, where $v^{\eps}$ on
$\overline{I_+^{\eps}(1)}$ is now  given by
\begin{equation}
v^{\eps}(x,t)=\biggl[c_r(\tau)\int\limits_{\om(\tau;0,0)}^L
b_s^{\eps}(\xi)S(\xi,\tau)
a_s^{\eps}(\om(0;\xi,\tau))\,d\xi\biggr]
\bigg|_{\tau=\theta(x,t)}S(x,t).
\end{equation}
Hence (\ref{eq:208}) in this case
 is true by (\ref{eq:eps}) and continuity of $\la$.
Therefore for $|S_5(t)|$  the
estimate (\ref{eq:209}) holds.

We now return to (\ref{eq:206}) and,
taking into account (\ref{eq:209}),  estimate
$|w^{\eps_1}-w^{\eps_2}|$
following the proof of   (\ref{eq:203}). As a result,
the bound (\ref{eq:*}) is true for sufficiently small~$\eps_2$.
\end{proof}

\begin{cla}
Provided $\eps_2$ is small enough,
(\ref{eq:*}) is true on $\overline{I_+^{\eps_1}(1)}$.
\end{cla}
\begin{proof}
Note that $w^{\eps_1}$ and $w^{\eps_2}$ on $\overline{I_+^{\eps_1}(1)}$
are defined by the same formula
(\ref{eq:34_1}). Therefore
\begin{equation}\label{eq:210}
\begin{array}{c}\displaystyle
\Bigl(w^{\eps_1}-w^{\eps_2}\Bigr)(x,t)=\Bigl[c_r(S_6^{\eps_1}-S_6^{\eps_2})\Bigr]
(\theta(x,t))+\Bigl[c_r(S_7^{\eps_1}-S_7^{\eps_2})\Bigr]
(\theta(x,t))+S_1(x,t).
\end{array}
\end{equation}
The
upper bound $C\eps_2$ for the absolute value of the
 second summand follows from
 (\ref{eq:61_1}).
To estimate the
first summand, we use the equality
$$
\Bigl(S_6^{\eps_1}-S_6^{\eps_2}\Bigr)
(t)
$$
$$
=\int\limits_{\om(t;0,t_1^*-\eps_1^-(\eps_2))}^L
\biggl[\Bigl(b_s^{\eps_1}+b_r\Bigr)w^{\eps_1}-
\Bigl(b_s^{\eps_2}+b_r\Bigr)w^{\eps_2}
+b_r\Bigl(v^{\eps_1}
-v^{\eps_2}\Bigr)\biggr](x,t)\,dx
$$
$$
+\int\limits_{\om(t;0,t_1^*-\eps_1^-(\eps_1))}^{\om(t;0,t_1^*-\eps_1^-
(\eps_2))}
\Bigl(b_rw^{\eps_1}\Bigr)(x,t)\,dx.
$$
The absolute value of the first
summand is already  estimated in
the proof of Claim 4. For the second summand we can apply Claim 2, since
$[\om(t;0,t_1^*-\eps_1^-(\eps_1)),
\om(t;0,t_1^*-\eps_1^-(\eps_2))]\times\{t\}$ $\subset$
$\overline{\Pi^{\eps_1}(1)}$
for $t\in[t_1^*-\eps_1^-(\eps_2),
t_1^*+\eps_1^+(\eps_2)]$.
As a consequence,
\begin{equation}\label{eq:S^1}
\Bigl|\Bigl(S_6^{\eps_1}-S_6^{\eps_2}\Bigr)(t)\Bigr|\le C\eps_2.
\end{equation}
Applying the method of sequential approximation to
$w^{\eps_1}-w^{\eps_2}$ given by   (\ref{eq:210}), on the account of
(\ref{eq:S^1}) and the upper bound $C\eps_2$ for the second summand
in (\ref{eq:210}), we derive the bound
\begin{equation}
\Bigl|\Bigl(w^{\eps_1}-w^{\eps_2}\Bigr)(x,t)\Bigr|\le C\eps_2
\exp\biggl\{T\Bigl(\max\limits_{(x,t,y)\in\overline{\Pi^T}\times\R}
|(\nabla_uf)(x,t,y)|
+\max\limits_{(x,t)\in\overline{\Pi^T}}|p(x,t)|\Bigr)\biggr\}.
\end{equation}
This implies (\ref{eq:*}), provided $\eps_2$ is chosen small enough.
\end{proof}

\begin{cla}
Provided $\eps_2$ is small enough,
(\ref{eq:*}) is true on $\overline{I_+^{\eps_2}(1)}\cap
\overline{\Pi^{\eps_1}(2)}$.
\end{cla}
\begin{proof}
We follow the proof of Claim 4 with the following changes.
On the account of (\ref{eq:supp}), for $w^{\eps_1}-w^{\eps_2}$
on $\overline{I_+^{\eps_2}(1)}\cap
\overline{\Pi^{\eps_1}(2)}$
we have the representation
(\ref{eq:206}) with $t_1^*-\eps_1^-(\eps_2)$  replaced by
$t_1^*+\eps_1^+(\eps_1)$ in the fifth summand,
with
\begin{equation}\label{eq:S_5}
S_5(t)=c_r(t)\biggl[\int\limits_{\om(t;0,t_1^*+\eps_1^+(\eps_1))}^
{\om(t;0,t_1^*-\eps_1^-(\eps_1))}
\Bigl(b_rv^{\eps_1}\Bigr)(x,t)\,dx-
\int\limits_0^{\om(t;0,t_1^*-\eps_1^-(\eps_2))}
\Bigl(
b_rv^{\eps_2}\Bigr)(x,t)\,dx\biggr],
\end{equation}
and with one more summand
\begin{equation}\label{eq:r8}
\biggl[ c_r(\tau)
\int\limits_{\om(\tau;0,t_1^*+
\eps_1^+(\eps_1))}^{\om(\tau;0,t_1^*-
\eps_1^-(\eps_2))}
\Bigl(b_r
(w^{\eps_1}-w^{\eps_2})\Bigr)(\xi,\tau)\,d\xi\biggr]\bigg|_{\tau=\theta(x,t)}.
\end{equation}
To estimate the absolute value of the latter expression, we use Claims 4
and 5.
To estimate $|S_5(t)|$, to both  integrals we apply the same
argument that was used for evaluation of $|S_5(t)|$ in  the proof of Claim 4.
\end{proof}

By Claims 3--6, Assertion 1 is true for $j=1$.

\begin{cla}
The functions $w^{\eps}$ are bounded
 on $\overline{I_+^{\eps}(1)}$, uniformly in $\eps>0$.
\end{cla}
\begin{proof}
The claim follows from (\ref{eq:34_1}), (\ref{eq:208}),
(\ref{eq:20}),
 and Assumptions 5 and 7.
\end{proof}

By Claims 2 and 7, Assertion 2 is true for $j=1$.

\begin{cla}
The family of functions $w^{\eps}$ converges in
$\Con(\overline{\Pi_0^T}\cup\overline{\Pi(1)})$ as $\eps\to 0$.
\end{cla}
\begin{proof}
Since, by Proposition 1, each $w^{\eps}$ is continuous, it suffices to prove
the convergence separately on $\overline{\Pi_0^T}$
 and $\overline{\Pi(1)}$. On the former domain the convergence is given by
Claim 2. The convergence on the latter domain follows by the Cauchy
criterion which holds by  Claims  3--5, and the fact that
$\overline{\Pi(1)}\subset
\overline{\Pi^{\eps_2}(1)}\cup
\overline{I_+^{\eps_2}(1)}$ for every
$\eps_2>0$.
Thus, Assertion 5 is true for $j=1$.

In the sequel the limit function will be denoted by $w(x,t)$.
\end{proof}

\begin{cla}
The estimate $D_1(\eps)\le C$
 is true for all $\eps>0$.
\end{cla}
\begin{proof}
{\it Case 1. $(0,\tilde t_1)\in I_+(1)$.} Then
$$
D_1(\eps)=\int\limits_{t_1^*-\eps_1^-(\eps)}
^{t_1^*+\eps_1^+(\eps)}\biggl|c_r(t)
\int\limits_{Q(t)}^Lb_s^{\eps}(x)S(x,t)a_s^{\eps}(\om(0;x,t))\,dx\biggr|\,dt
$$
$$
\le\max\limits_{(x,t)\in\overline{\Pi_0^T}}\biggl|
\frac{(c_rS)(x,t)}{(\d_t\om)(0;x,t)}
\biggr|\int\limits_0^L|a_s^{\eps}(x)|\,dx
\int\limits_0^L|b_s^{\eps}(x)|\,dx\le C.
$$
The estimate follows from
(\ref{eq:eps_2}) and
Assumptions 5 and 6.

{\it Case 2. $(0,t_1)\in I_+(1)$.} Then
$$
D_1(\eps)=\int\limits_{t_1-\eps}
^{t_1+\eps}|c_s^{\eps}(t)
(S_6^{\eps}+S_7^{\eps})(t)|\,dt\le C.
$$
This  estimate follows from  (\ref{eq:S_6}), (\ref{eq:61_1}),
and (\ref{eq:eps_2}).

Thus, Assertion 3 is true for $j=1$.
\end{proof}

\begin{cla}
If $\eps_2$ is small enough and $\eps_1<\eps_2$, then
$$
|R_1(\eps_1,\eps_2)|\le C\eps_2.
$$
\end{cla}
\begin{proof}
We will consider $\eps_2$ as small  as in Claims 3--6.
We will use (\ref{eq:9}) restricted to
$[t_1^*-\eps_1^-(\eps_2),t_1^*+\eps_1^+(\eps_2)]$ and
(\ref{eq:30_1}).

{\it Case 1. $(0,\tilde t_1)\in I_+(1)$.}
By (\ref{eq:r3}) and (\ref{eq:r4}) we have
$c_s^{\eps}=0$, and therefore
$$
R_1(\eps_1,\eps_2)=\int\limits_{t_1^*-\eps_1^-(\eps_2)}
^{t_1^*+\eps_1^+(\eps_2)}c_r(t)\int\limits_{Q(t)}^L
\Bigl(b_s^{\eps_1}v^{\eps_1}-b_s^{\eps_2}
v^{\eps_2}\Bigr)(x,t)\,dx\,dt.
$$
Applying (\ref{eq:30_1}) restricted to $\overline{\Pi_0^T}$ and
changing  coordinates $(x,t)\to (x,\xi)=(x,\om(0;x,t))$,
we obtain
$$
R_1(\eps_1,\eps_2)=\sum\limits_{(q,d)\in E}
\int\limits_{x_d^*-\eps_2}
^{x_d^*+\eps_2}\int\limits_{x_q-\eps_2}^{x_q+\eps_2}
Q(x,\xi)
\Bigl[b_s^{\eps_1}(x)a_s^{\eps_1}(\xi)-b_s^{\eps_2}(x)a_s^{\eps_2}(\xi)\Bigr]
\,dx\,d\xi,
$$
where
$$
Q(x,\xi)=\frac{(c_rS)(x,t)}{(\d_t\om)(0;x,t)}\Big|_{t=\tilde\om(x;\xi,0)}
$$
is a  continuous function in $x$ and $\xi$ and $E$
is the set of pairs of indices
$q\le k$ and $d\le m$ such that
$\om(0;x_q,\tilde t_1)=x_d^*$.
Evidently, if $(0,\tilde t_1)\in I_+(1)$, then there exists at least one pair
 $(q,d)$ that satisfies the latter condition.
Then
$$
|R_1(\eps_1,\eps_2)|\le\sum\limits_{q,d}\max\limits_{(x,\xi)\in[x_q-\eps_2,x_q+\eps_2]\times
[x_d^*-\eps_2,x_d^*+\eps_2]}|Q(x,\xi)-
Q(x_q,x_d^*)|
$$
$$
\times\biggl[\int\limits_{x_q-\eps_2}^{x_q+\eps_2}|b_s^{\eps_1}|\,dx
\int\limits_{x_d^*-\eps_2}^{x_d^*+\eps_2}|a_s^{\eps_1}|\,dx+
\int\limits_{x_q-\eps_2}^{x_q+\eps_2}|b_s^{\eps_2}|\,dx
\int\limits_{x_d^*-\eps_2}^{x_d^*+\eps_2}|a_s^{\eps_2}|\,dx\biggr]
$$
$$
+\sum\limits_{q,d}|Q(x_q,x_d^*)|\biggl|\int\limits_{x_q-\eps_2}^{x_q+\eps_2}b_s^{\eps_1}\,dx
\int\limits_{x_d^*-\eps_2}^{x_d^*+\eps_2}(a_s^{\eps_1}-
a_s^{\eps_2})\,dx+
\int\limits_{x_d^*-\eps_2}^{x_d^*+\eps_2}a_s^{\eps_2}\,dx
\int\limits_{x_q-\eps_2}^{x_q+\eps_2}(b_s^{\eps_1}-
b_s^{\eps_2})\,dx
\biggr|.
$$
The last two summands are equal to 0 by (\ref{eq:eps_1}).
The first summand is bounded from above by $C\eps_2$ by
(\ref{eq:eps_2}) and the continuity property for $Q(x,t)$. This
completes the proof in Case 1.

{\it Case 2. $(0,t_1)\in I_+(1)$.} Then the second summand in
(\ref{eq:9}) is equal to 0, and therefore
\begin{equation}\label{eq:R_i}
R_1(\eps_1,\eps_2)=\int\limits_{t_1-\eps_2}
^{t_1+\eps_2}\Bigl[c_s^{\eps_1}(t)
\Bigl(S_6^{\eps_1}+S_7^{\eps_1}\Bigr)(t)-c_s^{\eps_2}(t)
\Bigl(S_6^{\eps_2}+S_7^{\eps_2}\Bigr)
(t)\Bigr]\,dt.
\end{equation}
We hence have
\begin{equation}\label{eq:R}
\begin{array}{c}\displaystyle
|R_1(\eps_1,\eps_2)|=\biggl|\int\limits_{t_1-\eps_2}
^{t_1+\eps_2}\Bigl(c_s^{\eps_1}S_7^{\eps_1}-c_s^{\eps_2}S_7^{\eps_2}\Bigr)(t)
\,dt
+\int\limits_{t_1-\eps_2}
^{t_1+\eps_2}\Bigl(c_s^{\eps_1}-c_s^{\eps_2}\Bigr)(t)\Bigl(S_6^{\eps_1}(t)-
S_6^{\eps_1}(t_1)\Bigr)\,dt
\\
\displaystyle
+
S_6^{\eps_1}(t_1)\int\limits_{t_1-\eps_2}
^{t_1+\eps_2}\Bigl(c_s^{\eps_1}-c_s^{\eps_2}\Bigr)(t)\,dt
+\int\limits_{t_1-\eps_2}
^{t_1+\eps_2}c_s^{\eps_2}(t)\Bigl(S_6^{\eps_1}-S_6^{\eps_2}\Bigr)(t)
\,dt\biggr|
\\
\displaystyle
\le\int\limits_{t_1-\eps_2}
^{t_1+\eps_2}\Bigl(|c_s^{\eps_1}||S_7^{\eps_1}|+|c_s^{\eps_2}||S_7^{\eps_2}|
\Bigr)\,dt
+\int\limits_{t_1-\eps_2}
^{t_1+\eps_2}\Bigl(|c_s^{\eps_1}|+|c_s^{\eps_2}|
\Bigr)\Bigl|(S_6^{\eps_1}(t)-
S_6^{\eps_1}(t_1))\Bigr|
\,dt\\
\displaystyle
+\int\limits_{t_1-\eps_2}
^{t_1+\eps_2}|c_s^{\eps_2}|\Bigl|\Bigl(S_6^{\eps_1}-S_6^{\eps_2}\Bigr)(t)\Bigr|\,dt.
\end{array}
\end{equation}
Combining (\ref{eq:61_1}) with (\ref{eq:eps_2}),
we obtain the upper bound $C\eps_2$ for the first summand in the right-hand
side of the inequality (\ref{eq:R}). The same bound for the last summand
follows from
(\ref{eq:S^1}) and (\ref{eq:eps_2}).
It remains to estimate the second summand. Without loss of
generality we assume that $t\ge t_1$.
By (\ref{eq:supp}) we have
the representation
\begin{equation}\label{eq:76_i}
\begin{array}{c}\displaystyle
S_6^{\eps_1}(t)-S_6^{\eps_1}(t_1)=
\int\limits_{\om(t;0,t_1-\eps_1)}^L
\Bigl(
b_s^{\eps_1}+b_r\Bigr)(x)\Bigl(
w^{\eps_1}(x,t)-w^{\eps_1}(x,t_1)\Bigr)
\,dx\\
\displaystyle
+\int\limits_{\om(t_1;0,t_1-\eps_1)}^
{\om(t;0,t_1-\eps_1)}b_r(x)w^{\eps_1}(x,t_1)\,dx
+\int\limits_{Q(t)}^L(b_rS)(x,t)
a_s^{\eps_1}(\om(0;x,t))\,dx\\
\displaystyle
-\int\limits_{\om(t_1;0,0)}^
{\om(t_1;L,t)}\Bigl(b_rS\Bigr)(x,t_1)
a_s^{\eps_1}(\om(0;x,t_1))\,dx
-\int\limits_{\om(t_1;L,t)}^L\Bigl(b_rS\Bigr)(x,t_1)
a_s^{\eps_1}(\om(0;x,t_1))\,dx.
\end{array}
\end{equation}
Changing  coordinates
$x\to\om(0;x,t)$ in the third integral and
$x\to\om(0;x,t_1)$ in the fourth integral,
after simple computation we have
\begin{equation}\label{eq:77_i}
\begin{array}{c}\displaystyle
\Bigl|S_6^{\eps_1}(t)-S_6^{\eps_1}(t_1)\Bigr|
\le
\biggl(\int\limits_0^L|b_s^{\eps}|\,dx+L\max\limits_{x\in[0,L]}|b_r(x)|\biggr)
\\[8mm]
\displaystyle
\times\max\limits_{|t-t_1|\le\eps_1,(x,t)\in\overline{\Pi_0^T}\cup
\overline{\Pi(1)}}
\Bigl[\Bigl|\Bigl(w^{\eps_1}-w\Bigr)(x,t)\Bigr|+
\Bigl|w(x,t)-w(x,t_1)\Bigr|\\[8mm]\displaystyle
+\Bigl|w(x,t_1)-w^{\eps_1}(x,t_1)\Bigr|\Bigr]\\[8mm]
\displaystyle
+\max\limits_{(x,t)\in\overline{\Pi^{\eps_1}(1)}}|w^{\eps_1}(x,t)|
\max\limits_{x\in[0,L]}|b_r(x)|\Bigl|
\om(t;0,t_1-\eps_1)-\om(t_1;0,t_1-\eps_1)\Bigr|
\\[8mm]
\displaystyle
+\biggl[\max\limits_{(x,t)\in[0,L]\times[t_1-\eps_1,t_1+\eps_1]}\biggl|
\frac{(b_rS)(\om(t;x,0),t)}{(\d_x\om)(0;\om(t;x,0),t)}
-\frac{(b_rS)(\om(t_1;x,0),t_1)}{(\d_x\om)(0;\om(t_1;x,0),t_1)}\biggr|
\\[8mm]
\displaystyle
+\max\limits_{x\in[\om(t_1-\eps_1;L,t_1+\eps_1),L]}|b_r(x)|
\max\limits_{(x,t)\in\overline{\Pi^T}}
|S(x,t)|\biggr]\int\limits_0^L|a_s^{\eps_1}|\,dx\le C\eps_2.
\end{array}
\end{equation}
The latter estimate is true
by Claims 2, 7, and 8, estimates (\ref{eq:eps_2}), and Assumptions 3 and 5.

Claim 9 follows, and therefore
 Assertion 4 is true for $j=1$.
\end{proof}

{\bf Induction assumption.}
We assume that
Assertions 1--5 are true for  $j\le i-1$, $i\ge 2$.

\begin{cla}
The functions $w^{\eps}(x,t)$ are bounded on
$\overline{\Pi^{\eps}(i)}$, uniformly in $\eps>0$.
\end{cla}
\begin{proof}
The proof is similar to the proof of Claim 2. The function
$w^{\eps}(x,t)$ on $\overline{\Pi^{\eps}(i)}$ is defined by the formula
(\ref{eq:27_0}), where
$$
(Rw^{\eps})(x,t)=c_r(\theta(x,t))\biggl[\int\limits_0^{\om(\tau;0,
t_{i-1}^*+\eps_{i-1}^+(\eps))}\Bigl(b_s^{\eps}+b_r\Bigr)(\xi)
w^{\eps}(\xi,\tau)\,d\xi
$$
$$
+\int\limits_{\om(\tau;0,
t_{i-1}^*+\eps_{i-1}^+(\eps))}^L\Bigl(b_s^{\eps}+b_r\Bigr)(\xi)
w^{\eps}(\xi,\tau)\,d\xi+
\int\limits_{\om(\tau;0,
t_{i-1}^*+\eps_{i-1}^+(\eps))}^{Q(\tau)}(b_rS)(\xi,\tau)v^{\eps}
(0,\theta(\xi,\tau))\,d\xi
$$
$$
+\int\limits_{Q(\tau)}^L(b_rS)(\xi,\tau)a_s^{\eps}
(\om(0;\xi,\tau))\,d\xi\biggr]\bigg|_{\tau=\theta(x,t)}.
$$
Taking into account (\ref{eq:eps_2}) and
Assertions 2 and 3 for $j<i$, we conclude that the last three summands are
bounded uniformly in $\eps$.
Similarly to \cite[p.~646]{KmitHorm}, we obtain the global estimate
\begin{equation}
\begin{array}{c}\displaystyle
\max\limits_{(x,t)\in\overline{\Pi_{1}^{\eps}(i)}}|w^{\eps}(x,t)|\\
\displaystyle
\le\biggl(\frac{1}{1-q^1t^1}\biggr)^
{\lceil\frac{T}{t^1}\rceil} P(E)
\biggl(1+\max
\limits_{(x,t)\in\Bigl(\bigcup_{j=1}^{i-1}
\overline{\Pi(j)}\Bigr)\cup
\overline{I_+^{\eps}(i-1)}\cup
\overline{\Pi_0^T}}|w^{\eps}(x,t)|
\biggr),
\end{array}
\end{equation}
where $q^1$, $t^1$, and $E$
are defined
by  (\ref{eq:110}), (\ref{eq:t_0}), and (\ref{eq:48_1}), and
$P(E)$ is a
polynomial of degree
$\lceil T/t^1\rceil$
with  positive coefficients depending on $f(x,t,0),L$, and~$T$.
The claim now follows from Assertion 2 for $j\le i-1$.
\end{proof}


\begin{cla}
1. Provided $\eps_2$ is small enough,
for all $\eps'_2\le\eps_2$ and
for all $\eps'_1\le\eps'_2$ the estimate
(\ref{eq:r7})
is true on $\overline{\Pi(i)}
\setminus I_+^{\eps_2}(i)$.\\
2. Provided $\eps_2$ is small enough,
(\ref{eq:*}) is true on
$\overline{\Pi^{\eps_2}(i)}$.
\end{cla}
\begin{proof}
We fix an arbitrary sequence
$t_{i-1}^*=p_0<p_1<p_2<\dots<p_M=t_i^*-\eps_i^-(\eps_2)$
such that $M\ge 2$,
$p_1>t_{i-1}^*+\eps_{i-1}^+(\eps_2)$, and
 $\supp b_s^{\eps_2}\subset
[\om(p_i;0,p_{i-1}),L]$.
We can do so due to (\ref{eq:supp}).
Given this sequence,
we devide $\overline{\Pi(i)}
\setminus I_+^{\eps_2}(i)$
into a finite number of subsets
$$
\Pi^{\eps_2}(i,j)=\Bigl\{(x,t)\in
\overline{\Pi(i)}
\setminus I_+^{\eps_2}(i)\,|\,
\tilde\om(x;0,p_{j-1})\le t\le
\tilde\om(x;0,p_j)\Bigr\}.
$$
Note that, if $j<M$, then
$\Pi^{\eps_2}(i,j)$ actually does not depend
on $\eps_2$.

 To obtain (\ref{eq:r7}) for $\Pi^{\eps_2}(i,1)$, we first derive
(\ref{eq:*}) for $\Pi^{\eps_2}(i,1)$. Doing this, we follow the
proof of Claim 3 (Item 1) for Case 1
with the following
changes.   Throughout the proof $\eps'_1$ and $\eps'_2$ are
replaced by $\eps_1$ and $\eps_2$, respectively.
We use the formulas (\ref{eq:200})  and (\ref{eq:51_i})
with
$\om(t;0,t_{i-1}^*+\eps_{i-1}^+(\eps_2))$ in place of
$Q(t)$
and with
$$
S_2(t)=c_r(t)\int
\limits_{\om(t;0,t_{i-1}^*)}^L
\Bigl[(b_s^{\eps_1}-b_s^{\eps_2})
w^{\eps_2}
+b_s^{\eps_1}(w^{\eps_1}-w^{\eps_2})
\Bigr]
(x,t)\,dx
$$
$$
+c_r(t)\int\limits_{\om(t;0,
t_{i-1}^*+\eps_{i-1}^+(\eps_2))}^Lb_r(x)
\Bigl(w^{\eps_1}
-w^{\eps_2}\Bigr)
(x,t)\,dx.
$$
To estimate $|S_2|$, we use
Assertions 1, 2, and 5 for $j<i$.
To estimate the analog of
$|S_3|$, we use
Assertions 4 and 5 for $j<i$. As a result,
(\ref{eq:*}) is true on
$\Pi^{\eps_2}(i,1)$.

Note that $\Pi^{\eps_2}(i,1)$ does not depend on
$\eps_1$ and $\eps_2$. By Assertion 1 for $j=i-1$ and
(\ref{eq:*}) for $\Pi^{\eps_2}(i,1)$,
that we have just proved, we conclude that $w^{\eps}$
converges in $\Con(\Pi^{\eps_2}(i,1))$. The estimate
(\ref{eq:r7}) is hence true on $\Pi^{\eps_2}(i,1)$.

The estimate
(\ref{eq:r7}) for $\Pi^{\eps_2}(i,2)$ follows similarly to
 the proof of Claim 3 (Item 1) for Case 2 (see the proof of
(\ref{eq:r7}) for $\Pi^{\eps_2}(1,2)$). The minor change  is that
we estimate $|S_2(t)+S_4(t)|$ using Assertion 5 for $j=i-1$ and
(\ref{eq:r7}) for $\Pi^{\eps_2}(i,1)$, that we have just proved.
Similar arguments apply to the subsets
 $\Pi^{\eps_2}(i,j)$, $j\le M$. As a result, (\ref{eq:r7}) is true
for the whole $\overline{\Pi(i)}
\setminus I_+^{\eps_2}(i)$.

Item 1 follows.

Item 2 is straightforward consequence of Item 1.
\end{proof}

\begin{cla}
Provided $\eps_2$ is small enough,
(\ref{eq:*}) is true on
$\overline{\Pi^{\eps_1}(i)}\cap\overline{I_+^{\eps_2}(i)}$.
\end{cla}
\begin{proof}
The proof is similar to the proof of Claim 4. We concentrate
only on the
changes that here appear.
We choose $\eps_2$ so small  that the condition (\ref{eq:supp})
with $t_1^*+\eps_1^+(\eps_2)$  replaced by $t_i^*+\eps_i^+(\eps_2)$
and  $t_1^*-\eps_1^-(\eps_2)$  replaced by $t_i^*-\eps_i^-(\eps_2)$
is true.
As in Claim 4, we distinguish two cases.

{\it Case 1. There exists $j\le l$ such that $(0,t_j)\in I_+(i)$}.
We use (\ref{eq:206}) with
$t_{j}^*$ in place of
$t_1^*$, where $Q(t)$ in $S_3(t)$  is replaced by
$\om(t;0,t_{j}-\eps_2)$.
For the
first four  summands in the analog of (\ref{eq:206}) we use
Assertions 1--4 for $j<i$ and
 Claim 12.
The upper bound for $|S_5(t)|$ follows from the bounds
\begin{equation}\label{eq:64_i}
\begin{array}{c}\displaystyle
|S_6^{\eps}(t)|\le C,\quad t\in[t_i^*-\eps_i^-(\eps),t_i^*+\eps_i^+(\eps)]
\end{array}
\end{equation}
and
\begin{equation}\label{eq:68_i}
\begin{array}{c}\displaystyle
|S_7^{\eps}(t)|\le C\eps,\quad t\in[t_i^*-\eps_i^-(\eps),t_i^*+\eps_i^+(\eps)],
\end{array}
\end{equation}
where $S_6^{\eps}(t)$ and $S_7^{\eps}(t)$ are defined by (\ref{eq:SS})
with $t_i^*-\eps_i^-(\eps)$ in place of $t_1^*-\eps_1^-(\eps)$ and
$\om(t;0,t_{i-1}^*+\eps_{i-1}^+(\eps))$ in place of $Q(t)$.
The estimate (\ref{eq:64_i}) is true by the arguments used for obtaining
(\ref{eq:S_6}), Assertions 2 and 3 for $j<i$, and Claim 11.
The estimate (\ref{eq:68_i}) is obtained similarly to
(\ref{eq:61_1}). The estimates (\ref{eq:64_i}) and (\ref{eq:68_i})
imply (\ref{eq:208}).

{\it Case 2.    $(0,t_j)\not\in I_+(i)$ for $j\le l$}.
On the account of (\ref{eq:supp}), (\ref{eq:r3}), and (\ref{eq:r4}), on
$\overline{I_+^{\eps}(i)}$
we have the equality
\begin{equation}\label{eq:V}
\begin{array}{c}\displaystyle
v^{\eps}(x,t)=
\biggl[c_r(\tau)\int\limits
_{\om(\tau;0,t_{i-1}^*+\eps_{i-1}^+(\eps_2))}^{Q(\tau)}
b_s^{\eps}(\xi)S(\xi,\tau)
v^{\eps}
(0,\theta(\xi,\tau))\,d\xi\biggr]\bigg|_
{\tau=\theta(x,t)}S(x,t)
\\
\displaystyle
+\biggl[c_r(\tau)\int\limits_{Q(\tau)}^L
b_s^{\eps}(\xi)S(\xi,\tau)
a_s^{\eps}
(\om(0;\xi,\tau))\,d\xi\biggr]\bigg|_
{\tau=\theta(x,t)}S(x,t).
\end{array}
\end{equation}
For the second summand we apply the same arguments as in the proof of Claim 4
for Case 2.
For the first summand we apply (\ref{eq:eps}), Assertion 3 for $j<i$ and
the inclusion
$$
\{(0,\theta(x,t))\,|\,t\in[t_i-\eps_i^-(\eps),t_i+\eps_i^+(\eps)],
x\in[\om(t;0,t_{i-1}^*+\eps_{i-1}^+(\eps_2)),Q(t)]\}
$$
$$
\subset\bigcup\limits_{j=1}^{i-1}\Bigl(
\overline{I_+^{\eps}(j)}\cap\{(x,t)\,|\,x=0\}\Bigr).
$$
Therefore
(\ref{eq:208}) is true on $\overline{I_+^{\eps}(i)}$ . The claim follows.
\end{proof}

From what has already been proved we conclude that (\ref{eq:208}) holds on
$\overline{I_+^{\eps}(j)}$ for every $j\le n(T)$.
Since $v^{\eps}(x,t)=S(x,t)a_s^{\eps}(\om(0;x,t))$ for $(x,t)\in
\overline{\Pi_0^T}$, (\ref{eq:208}) holds on $\overline{\Pi_0^T}$.

By (\ref{eq:0}),
$\supp v^{\eps}\subset\overline I_+^{\eps}$. Therefore
\begin{equation}\label{eq:1eps}
v^{\eps}=O\biggl(\frac{1}{\eps}\biggr)\quad \mbox{on}\quad
\overline{\Pi^T}.
\end{equation}
We will need this property in the sequel.

\begin{cla}
Provided $\eps_2$ is small enough,
(\ref{eq:*}) is true on
$\overline{I_+^{\eps_1}(i)}$.
\end{cla}
\begin{proof}
We follow the proof  of Claim 5 with the following
changes. We use (\ref{eq:SS}) and (\ref{eq:210})
with $t_1^*-\eps_1^-(\eps)$ and $Q(t)$
to be replaced
by $t_i^*-\eps_i^-(\eps)$.
To estimate the absolute value of the first summand in the
analog of (\ref{eq:210}),
we use
Assertions 1--3 for $j<i$ and
Claims 11--13.
As a consequence,
\begin{equation}\label{eq:74_i}
\begin{array}{c}\displaystyle
\Bigl|\Bigl(S_6^{\eps_1}-S_6^{\eps_2}\Bigr)(x,t)\Bigr|\le C\eps_2,
\quad t\in[t_i^*-\eps_i^-(\eps),
t_i^*+\eps_i^+(\eps)].
\end{array}
\end{equation}
To estimate the second summand, we use (\ref{eq:68_i}).
\end{proof}

\begin{cla}
Provided $\eps_2$ is small enough,
(\ref{eq:*}) is true on
$\overline{I_+^{\eps_2}(i)}\cap\overline{\Pi^{\eps_1}(i+1)}$.
\end{cla}
\begin{proof}
We follow the proof of Claim 6 with the following changes.
For $w^{\eps_1}-w^{\eps_2}$ on
$\overline{I_+^{\eps_2}(i)}\cap\overline{\Pi^{\eps_1}(i+1)}$
we use
(\ref{eq:206}) with
$t_1^*-\eps_1^-(\eps_2)$ replaced by $t_i^*-\eps_i^-(\eps_2)$
in the fifth summand
and with one more summand (\ref{eq:r8}), where
$t_1^*+\eps_1^+(\eps_1)$ and
$t_1^*-\eps_1^-(\eps_2)$ are replaced by
$t_i^*+\eps_i^+(\eps_1)$ and
$t_i^*-\eps_i^-(\eps_2)$, respectively.
In the representation of $S_3(t)$ we
now have
$\om(t;t_{i-1}^*-\eps_{i-1}^-(\eps_2))$ in place of $Q(t)$.
For $S_5(t)$ we now use the formula (\ref{eq:S_5}) with
$t_i^*-\eps_i^-(\eps_1)$,
 $t_i^*+\eps_i^+(\eps_1)$,
and $t_i^*-\eps_i^-(\eps_2)$
in place of
$t_1^*-\eps_1^-(\eps_1)$,
 $t_1^*+\eps_1^+(\eps_1)$,
and $t_1^*-\eps_1^-(\eps_2)$,
respectively.

To estimate the absolute value of the analog of (\ref{eq:r8})
 we use Claims 13
and 14.
\end{proof}
By Claims 12--15, Assertion 1 is true for $j=i$.

\begin{cla}
The functions $w^{\eps}(x,t)$ are bounded on
$\overline{I_+^{\eps}(i)}$, uniformly in $\eps>0$.
\end{cla}
\begin{proof}
Note that the function $w^{\eps}(x,t)$ on
$\overline{I_+^{\eps}(i)}$    is defined by (\ref{eq:34_1})
and (\ref{eq:SS})
with $t_i^*-\eps_i^-(\eps)$ in place of $t_1^*-\eps_1^-(\eps)$ and
$\om(t;0,t_{i-1}^*+\eps_{i-1}^+(\eps))$ in place of $Q(t)$.
The  claim follows from (\ref{eq:64_i}), (\ref{eq:68_i}),
and Assumptions 5 and 7.
\end{proof}
Assertion 2 for $j=i$ follows from Claims 11 and 16.

\begin{cla}
The family of functions $w^{\eps}$ converges in
$\Con(\bigcup_{j=1}^i\overline{\Pi(j)}\cup\overline{\Pi_0^T})$ as $\eps\to 0$.
\end{cla}
\begin{proof}
Since, by Proposition 1, each $w^{\eps}$ is continuous, it suffices to prove
the convergence separately on $\overline{\Pi_0^T}$
 and $\bigcup_{j=1}^i\overline{\Pi(j)}$.
On the former domain the convergence is ensured by
Claim 2. The convergence on the latter domain follows by the Cauchy
criterion which holds by  Assertion  1 for $j\le i$,  and the fact that
$\bigcup_{j=1}^i\overline{\Pi(j)}
\subset\bigcup_{j=1}^i(\overline{\Pi^{\eps_2}(j)}\cup I_+^{\eps_2}(j))$
for every $\eps_2>0$.
\end{proof}
By Claim 17, Assertion 5 holds for $j=i$.

\begin{cla}
The estimate $D_i(\eps)\le C$
 is true for all $\eps>0$.
\end{cla}
\begin{proof}
{\it Case 1. $(0,\tilde t_i)\in I_+(i)$.} We have
$$
D_i(\eps)
=\int\limits_{t_i^*-\eps_i^-(\eps)}
^{t_i^*+\eps_i^+(\eps)}\biggl|c_r(t)
\int\limits_0^{Q(t)}b_s^{\eps}(x)S(x,t)
v^{\eps}(0,\theta(x,t))\,dx
$$
$$
+c_r(t)
\int\limits_{Q(t)}^Lb_s^{\eps}(x)S(x,t)a_s^{\eps}(\om(0;x,t))\,dx\biggr|\,dt
$$
$$
\le\max\limits_{(x,t)\in\overline{\Pi_1^T}}\biggl|
\frac{(c_rS)(x,t)}{(\d_t\theta)(x,t)}
\biggr|\sum\limits_{j=1}^{i-1}
\int\limits_{t_j^*-\eps_j^-(\eps)}
^{t_j^*+\eps_j^+(\eps)}|v^{\eps}(0,t)|\,dt
\int\limits_0^L|b_s^{\eps}(x)|\,dx
$$
$$
+\max\limits_{(x,t)\in\overline{\Pi_0^T}}\biggl|
\frac{(c_rS)(x,t)}{(\d_t\om)(0;x,t)}
\biggr|\int\limits_0^L|a_s^{\eps}(x)|\,dx
\int\limits_0^L|b_s^{\eps}(x)|\,dx\le C,
$$
where $Q(t)$ is defined by (\ref{eq:Q}).
This estimate is true by
(\ref{eq:eps_2}) and Assertion 3 for $j<i$.

{\it Case 2. $(0,t_j)\in I_+(i)$ for some $j\le l$.} We have
$$
D_i(\eps)=\int\limits_{t_j-\eps}
^{t_j+\eps}|c_s^{\eps}(t)
(S_6^{\eps}+S_7^{\eps})(t)|\,dt\le C,
$$
This  estimate is true by
 (\ref{eq:64_i}),
(\ref{eq:68_i}),
and (\ref{eq:eps_2}).
\end{proof}
Claim 18 implies
Assertion 3 for $j=i$.

\begin{cla}
The estimate $|R_i(\eps_1,\eps_2)|\le C\eps_2$
is true for $\eps_2$ so small
that
Assertion 1 holds for all $j\le i$.
\end{cla}
\begin{proof}
We follow the proof of Claim 10 with
the changes listed below.
Similarly to Claim 10, we distinguish two cases.

{\it Case 1. $(0,\tilde t_i)\in I_+(i)$.} We have
$R_i(\eps_1,\eps_2)=R_i^1(\eps_1,\eps_2)+R_i^2(\eps_1,\eps_2)$, where
$$
R_i^1(\eps_1,\eps_2)=\int\limits_{t_i^*-\eps_i^-(\eps_2)}
^{t_i^*+\eps_i^+(\eps_2)}c_r(t)\int\limits_{0}^{Q(t)}
\Bigl(b_s^{\eps_1}v^{\eps_1}
-b_s^{\eps_2}
v^{\eps_2}\Bigr)(x,t)\,dx\,dt
$$
and
$$
R_i^2(\eps_1,\eps_2)=\int\limits_{t_i^*-\eps_i^-(\eps_2)}
^{t_i^*+\eps_i^+(\eps_2)}c_r(t)\int\limits_{Q(t)}^L
\Bigl(b_s^{\eps_1}v^{\eps_1}-b_s^{\eps_2}
v^{\eps_2}\Bigr)(x,t)\,dx\,dt.
$$
If $Q(t)<L$, then
$[Q(t),L]\times\{t\}$ $\subset$
$\overline{\Pi_0^T}$. Therefore
 $|R_i^2(\eps_1,\eps_2)|$ can be estimated in the same way as  $|R_1(\eps_1,\eps_2)|$
was estimated in the
proof of Claim 10 for Case 1. The minor change is that now $E$ will
denote the set of pairs of indices
$q\le k$ and $d\le m$ such that
$\om(0;x_q,\tilde t_i)=x_d^*$. It remains to estimate $|R_i^1(\eps_1,\eps_2)|$.
Applying  (\ref{eq:30}) and (\ref{eq:30_1}) restricted
to $\overline{\Pi_1^T}$ and
changing  coordinates $(x,t)\to (x,\xi)=(x,\theta(x,t))$,we obtain
$$
R_i^1(\eps_1,\eps_2)=\sum\limits_{(q,d)\in J}
\int\limits_{t_d^*-\eps_d^-(\eps_2)}
^{x_d+\eps_d^+(\eps_2)}\int
\limits_{x_q-\eps_2}^{x_q+\eps_2}
Q_1(x,t)
\Bigl[(b_s^{\eps_1}-b_s^{\eps_2})(x)v^{\eps_1}(0,t)
$$
$$
-b_s^{\eps_2}(x)
(v^{\eps_1}-v^{\eps_2})(0,t)\Bigr]
\,dx\,dt,
$$
where
$$
Q_1(x,\xi)=\frac{(c_rS)(x,\tau)}
{(\d_t\theta)(x,\tau)}\Big|_{\tau=\tilde\om(x;0,t)}
$$
and $J$ is the set of pairs of indices
$q\le k$ and $d\le i-1$ such that
$\om(x_q;0,t_d^*)=\tilde t_i$. Obviously,
at least one of the sets $E$ or $J$ is nonempty.
To estimate $|R_i^1(\eps_1,\eps_2)|$, in addition to the arguments used
for estimation of
$|R_2(\eps_1,\eps_2)|$
we apply Assertions 3 and 4 for $j<i$.

{\it Case 2. $(0,t_j)\in I_+(i)$ for some $j\le l$.}
We use
(\ref{eq:R_i}) and (\ref{eq:R}) with $t_1$ replaced by $t_j$. To estimate
the first and the third summands in the analog of (\ref{eq:R}), we apply
(\ref{eq:eps_2}), (\ref{eq:68_i}), and (\ref{eq:74_i}).
To estimate the second
summand, we use the representation
$$
S_6^{\eps_1}(t)-S_6^{\eps_1}(t_j)=
\int\limits_{\om(t;0,t_j-\eps_1)}^L
\Bigl(
b_s^{\eps_1}+b_r\Bigr)(x)\Bigl(
w^{\eps_1}(x,t)-w^{\eps_1}(x,t_j)\Bigr)
\,dx
$$
$$
+\int\limits_{\om(t_j;0,t_j-\eps_1)}^
{\om(t;0,t_j-\eps_1)}b_r(x)w^{\eps_1}(x,t_j)\,dx
$$
$$
+\biggl[\int\limits_{Q(t)}^L(b_rS)(x,t)
a_s^{\eps_1}(\om(0;x,t))\,dx-
\int\limits_{Q(t_j)}^L(b_rS)(x,t_j)
a_s^{\eps_1}(\om(0;x,t_j))\,dx\biggr]
$$
\begin{equation}
\begin{array}{c}\displaystyle
+\biggl[\int\limits_{\om(t;0,t_j-\eps_1)}^{Q(t)}(b_rS)(x,t)
v^{\eps_1}(0,\theta(x,t))\,dx-
\int\limits_{\om(t_j;0,t_j-\eps_1)}^{Q(t_j)}(b_rS)(x,t_j)
v^{\eps_1}(0,\theta(x,t_j))\,dx\biggr].
\end{array}
\end{equation}
The absolute value of the first three
summands are estimated similarly to estimation of
$|S_6^{\eps_1}(t)-S_6^{\eps_1}(t_1)|$ in the proof of Claim 10
(see (\ref{eq:77_i})).
In addition to the arguments used for (\ref{eq:77_i}), we apply
Assertion 2 for $j<i$ and Claims 11, 16, and 17.
We now concentrate on the last summand.
Let us rewrite it in the form
$$
\int\limits_{\theta(Q(t),t)}
^{t_j-\eps_1}\biggl[\frac{(b_rS)(x,t)}
{(\d_x\theta)(x,t)}\bigg|_{x=\om(t;0,\tau)}-
\frac{(b_rS)(x,t_j)}
{(\d_x\theta)(x,t_j)}\bigg|_{x=\om(t_j;0,\tau)}\biggr]
v^{\eps_1}(0,\tau)\,d\tau
$$
$$
+\int\limits_{\theta(Q(t_j),t_j)}
^{\theta(Q(t),t)}\frac{(b_rS)(x,t)}
{(\d_x\theta)(x,t)}\bigg|_{x=\om(t;0,\tau)}
v^{\eps_1}(0,\tau)\,d\tau.
$$
The absolute value of this expression is less than or equal to
$$
\biggl[\max\limits_{(x,t)\in[0,L]\times[t_j-\eps_1,t_j+\eps_1]}
\biggl|\frac{(b_rS)(x,t)}
{(\d_x\theta)(x,t)}\bigg|_{x=\om(t;0,\tau)}-
\frac{(b_rS)(x,t_j)}
{(\d_x\theta)(x,t_j)}\bigg|_{x=\om(t_j;0,\tau)}\biggr|
$$
$$
+\max\limits_{x\in[\om(t_j-\eps_1;Q(t_j+\eps_1),t_j+\eps_1),Q(t_j-\eps_1)]}
|b_r(x)|\max\limits_{(x,t)\in\overline{\Pi^T}}|S(x,t)|\biggr]
$$
$$
\times\sum\limits_{r=1}^{i-1}\int\limits_{t_r^*-\eps_r^-(\eps_1)}^
{t_r^*+\eps_r^+(\eps_1)}|v^{\eps_1}
(0,t)|\,dt\le C\eps_2.
$$
The latter bound is true by Assertion 3 for $j<i$ and Assumption 3.

We conclude that Assertion 4 holds for $j=i$.
\end{proof}

Thus the induction step is done and the proof of Lemma 3 is complete.

\section{Proof of Lemma 4}

Given an arbitrary test function $\varphi\in\D(\Pi)$, we have
to show the convergence of
$<v^{\eps}(x,t),\varphi(x,t)>$ as $\eps\to 0$.
 Fix $T=T(\vphi)>0$ such that
$\supp\vphi\subset\Pi^T$.
Let
\begin{equation}\label{eq:infty}
v^{\eps}(x,t)=\sum\limits_{i=0}^{\infty}v_i^{\eps}(x,t),
\end{equation}
where
$$
\supp v_0^{\eps}(x,t)=
\supp\Bigl\{v^{\eps}(x,t)
\Big|_{\overline{\Pi_0}\cap
\overline{I_+^{\eps}}}\Bigr\},\quad
\Pi_0=\{(x,t)\in\Pi\,|\,x>\om(t;0,0)\},
$$
and
$$
\supp\Bigl\{v_i^{\eps}(x,t)
\Big|_{\overline{\Pi^T}}\Bigr\}=
\supp\Bigl\{v^{\eps}(x,t)
\Big|_{\overline{I_+^{\eps}(i)}}
\Bigr\}\quad
\mbox{for}\quad i\le \rho(T),
$$
$\rho(T)$ is as  in Definition~3.
 Clearly, $v_i^{\eps}(x,t)$
for $i\ge 1$
is supported  on one of the connected components of
$(\overline{\Pi}\setminus\Pi_0)\cap\overline{I_+^{\eps}}$.
The representation (\ref{eq:infty}) is true because
$\supp v^{\eps}\subset\overline{I_+^{\eps}}$ by (\ref{eq:0}).
Since $<v^{\eps}(x,t),\varphi(x,t)>=
<\sum\limits_{i=0}^{\rho(T)}v_i^{\eps}(x,t),\varphi(x,t)>$,
  it
suffices to prove that
$<v_i^{\eps}(x,t),\varphi(x,t)>$
converges as $\eps\to 0$ separately for each $0\le i\le \rho(T)$.

\begin{clai}
$<v_0^{\eps}(x,t),\varphi(x,t)>$ converges as $\eps\to 0$.
\end{clai}
\begin{proof}
By (\ref{eq:30}) and (\ref{eq:30_1}),
\begin{equation}\label{eq:32}
v_0^{\eps}(x,t)=S(x,t)a_{s}^{\eps}(\om(0;x,t)).
\end{equation}
Let us compute the action
$$
<v_0^{\eps}(x,t),\varphi(x,t)>=\int\limits_{\Pi_0^T} S(x,t)a_s^{\eps}(\om(0;x,t))\varphi(x,t)\,d(x,t)
$$
$$
=\int\limits_{\tilde\Pi_0^T}
\frac{S(\xi,t)
\varphi(\xi,t)}{(\d_x\om)(0;\xi,t)}\Big|_{\xi=\om(t;x,0)}a_s^{\eps}(x)
\,d(x,t),
$$
where
$$
\tilde\Pi_0^T=\Bigl\{(x,t)\in\R^2\,|\,
(\om(t;x,0),t)\in\Pi_0^T\Bigr\}.
$$
It follows easily that
$$
<v_0^{\eps}(x,t),\varphi(x,t)>
\longrightarrow_{\eps\to 0}\sum\limits_{i=1}^{m}\int
\limits_{\tilde\Pi_0^T\cap\{(x,t)\,|\,x=x_i^*\}} \frac{S(\xi,t)
\varphi(\xi,t)}{\om_{\xi}(0;\xi,t)}\Big|_{\xi=\om(t;x_i^*,0)}\,dt,
$$
Here we used a simple change of coordinates $F_0:(x,t)\to(\om(0;x,t),t)$,
where
 $F_0$ maps $\Pi_0^T$ to $\tilde\Pi_0^T$. Since $\tilde\Pi_0^T$ is a bounded
domain,
the claim is proved.
\end{proof}

\begin{clai}
$<v_1^{\eps}(x,t),\varphi(x,t)>$ converges as $\eps\to 0$.
\end{clai}
\begin{proof}
Two cases are possible.

{\it Case 1.    $(0,t_1)\in I_+(1)$}.
Let
$$
S_8^{\eps}(t)=\int\limits_{0}^{L}[(b_s^{\eps}+b_r)w^{\eps}+b_r
v^{\eps}]\,dx.
$$
From  (\ref{eq:6}) we conclude that
$c_r(t)S_8^{\eps}(t)=w^{\eps}(0,t)$.
By Lemma 3, the family of functions $w^{\eps}(0,t)$ is uniformly convergent on $[0,T]$.
Since $c_r(t)$ is an arbitrary continuous function, the same
assertion is
true for the family of functions
$S_8^{\eps}(t)$.
Hence there exists a continuous
function $S_8^0(t)$ such that
$$
\lim\limits_{\eps\to 0}S_8^{\eps}(t)=
S_8^{0}(t) \,\,\mbox{in}\,\, \Con[0,T].
$$
Therefore, if $|t-t_1|\le C\eps_2$, we have
$$
\Bigl|S_8^{\eps}(t)-
S_8^0(t_1)\Bigr|
\le
\Bigl|\Bigl(S_8^{\eps}-S_8^0\Bigr)
(t)\Bigr|
+\Bigl|S_8^0(t)
-S_8^0
(t_1)\Bigr|
\le C\eps_2.
$$
Note that $S_8^{\eps}(t)=S_6^{\eps}(t)+S_7^{\eps}(t)$
whenever $t\in[t_1-\eps,t_1+\eps]$.
Using
the representation
(\ref{eq:205}), we conclude that
$$
<v_1^{\eps}(x,t),\vphi(x,t)>
=
\int\limits_{I_+^{\eps_0}(1)}
c_s^{\eps}(\theta(x,t))
S_8^{\eps}
(\theta(x,t))S(x,t)\vphi(x,t)\,d(x,t)
$$
$$
=\int\limits_{\tilde I_+^{\eps_0}(1)}
c_s^{\eps}(t)
S_8^{\eps}(t)
\frac{S(x,\tau)\vphi(x,\tau)}
{(\d_t\tilde\om)(0;x,\tau)}\bigg|_{\tau=\tilde\om(x;0,t)}\,d(x,t)
$$
$$
\longrightarrow_{\eps\to 0}S_8^0(t_1)
\int\limits_{\tilde I_+^{\eps_0}(1)\cap\{(x,t)\,|\,t=t_1\}}
\frac{S(x,\tau)\vphi(x,\tau)}
{(\d_t\tilde\om)(0;x,\tau)}\bigg|_{\tau=\tilde\om(x;0,t_1)}
\,dx.
$$
Here
$$
\tilde I_+^{\eps_0}(1)=\Bigl\{(x,t)\in\R^2\,|\,(x,\tilde\om(x;0,t))
\in I_+^{\eps_0}(1)\Bigr\}
$$
is a bounded domain. This implies that the latter integral is finite, and therefore
the claim in this case is true.

{\it Case 2.    $(0,t_1)\not\in I_+(1)$}. Using the representation
of $v_1^{\eps}$ given by  (\ref{eq:V}),
consider the action
$$
<v_1^{\eps}(x,t),\vphi(x,t)>
$$
$$
=
<c_r(\theta(x,t))\int\limits_{\om(\theta(x,t);0,0)}^Lb_s^{\eps}(\xi)
a_s^{\eps}(\om(0;\xi,\theta(x,t)))
S(\xi,\theta(x,t))\,d\xi\, S(x,t),\vphi(x,t)>
$$
$$
=\int\limits_{I_+^{\eps_0}(1)}
\int\limits_{\om(\theta(x,t);0,0)}^L c_r(\theta(x,t))b_s^{\eps}(\xi)a_s^{\eps}
(\om(0;\xi,\theta(x,t)))S(\xi,\theta(x,t))S(x,t)\vphi(x,t)\,d\xi\,d(x,t).
$$
Changing coordinates
$
(\xi,x,t)\to (\xi,x,
\om(0;\xi,\theta(x,t))),
$
we convert  the latter expression into
$$
\int\limits_{\Omega}b_s^{\eps}(\xi)
a_s^{\eps}(t)
\biggl[c_r(\theta(x,\tau))S(\xi,\theta(x,\tau))S(x,\tau)\vphi(x,\tau)
$$
$$
\times
\Bigl[(\d_t\om)(0;\xi,\theta(x,\tau))
(\d_t\theta)(x,\tau)\Bigr]^{-1}\biggr]
\bigg|_{\tau=\tilde\om(x;0,\tilde\om(\xi;t,0))}
\,d(\xi,x,t),
$$
where
$$
\Omega=\Bigl\{(\xi,x,t)\,|\,t=\om(0;\xi,\theta(x,\tau)),
(x,\tau)\in I_+^{\eps_0}(1),
\om(\theta(x,\tau);0,0)<\xi<L\Bigr\}
$$
is a bounded domain, $I(i)$ is the set of those pairs  $(j,r)$
that $j\le k$, $r\le m$, and
 $\om(0;x_j,t_i^*)=x_r^*$.
It is now clear that, as $\eps\to 0$,
$<v_1^{\eps}(x,t),\vphi(x,t)>$ converges to
$$
c_r(t_1^*)\sum\limits_{(j,r)\in I(1)}
\frac{S(x_j,t_1^*)}{(\d_t\om)(0;x_j,t_1^*)}
\int\limits_{\Omega\cap\{(\xi,x,t)\,|\,\xi=x_j,t=x_r^*\}}
\frac{S(x,t)
\vphi(x,t)}
{(\d_t\theta)(x,t)}\bigg|_{t=\tilde\om(x;0,t_1^*)}\,dx.
$$
Since the  sum  is finite, the claim
follows.
\end{proof}

\begin{clai}
$\lim\limits_{\eps\to 0}v_1^{\eps}(0,t)=
C_1\de(t-t_1^*)$ in $\D'(\R_+)$,
where  $C_1$ is a real constant.
\end{clai}
\begin{proof}
Take a test function $\psi(t)\in\D(\R_+)$ and compute the
action $<v_1^{\eps}(0,t),\psi(t)>$. Similarly to the proof of Claim 2,
we  obtain that
$$
v_1^{\eps}(0,t)\longrightarrow_{\eps\to 0}
S_8^0(t_1)S(0,t_1)\de(t-t_1),\quad \mbox{if}\quad
(0,t_1)\in I_+(1),
$$
and
$$
v_1^{\eps}(0,t)\longrightarrow_{\eps\to 0}c_r(t_1^*)
\sum\limits_{\{j\le k\,|\,(j,r)\in I(1),r\le m\}}\frac{S(x_j,t_1^*)}{(\d_t\om)(0;x_j,t_1^*)}
S(0,t_1^*)\de(t-t_1^*),
$$
$$
 \mbox{if}\quad
(0,t_1)\not\in I_+(1).
$$
\end{proof}

\begin{clai}
$<v_j^{\eps}(x,t),\varphi(x,t)>$ for
$1\le j\le \rho(T)$ converges as $\eps\to 0$.
\end{clai}
\begin{proof}
We prove the claim, using induction on $j$.
The base case of $j=1$ is given
by Claims 2 and 3.
We make the following assumptions
 for $j\le i-1$.\\
{\it Assumption 1.}
$<v_j^{\eps}(x,t),\varphi(x,t)>$
converges as $\eps\to 0$.\\
{\it Assumption 2.}
$\lim\limits_{\eps\to 0}v_j^{\eps}(0,t)=
C_j\de(t-t_j^*)$ in $\D'(\R_+)$,
where  $C_j$ is a real constant.

 Prove the claim for
$j=i$.
Similarly to Claim 2, we distinguish
two cases.

{\it Case 1.    $(0,t_j)\in I_+(i)$ for some $j\le l$}.
The claim follows similarly to
the proof of Claim 2 for Case 1.

{\it Case 2.    $(0,t_j)\not\in I_+(i)$,  for $j\le l$}. We
use  (\ref{eq:V}),
 (\ref{eq:30}),
and (\ref{eq:30_1}), and represent $v_i^{\eps}$ in the form
\begin{equation}
\begin{array}{c}
\displaystyle
v_i^{\eps}(x,t)=\biggl[c_r(\tau)
\int\limits_
{\om(\tau;0,t_{i-1}^*+\eps_{i-1}^+
(\eps_0))}^{Q(\tau)}
b_s^{\eps}(\xi)v^{\eps}(0,\theta(\xi,\tau))S(\xi,\tau)
\,d\xi\biggr]
\bigg|_{\tau=\theta(x,t)}S(x,t)\nonumber\\
\displaystyle
+\biggl[c_r(\tau)
\int\limits_{Q(\tau)}^L
\Bigl(b_s^{\eps}v_0^{\eps}\Bigr)(\xi,\tau)
\,d\xi\biggr]
\bigg|_{\tau=\theta(x,t)}S(x,t),
\nonumber
\end{array}
\end{equation}
where $Q(t)$ is defined by
(\ref{eq:Q}).
The convergence of the second summand
 follows from Claim 1. We now prove the
convergence of the first summand. Consider the action
$$
<c_r(\theta(x,t))\int\limits_{
\om(\theta(x,t);0,t_{i-1}^*+\eps_{i-1}^+(\eps))}^{Q(\theta(x,t))}
b_s^{\eps}(\xi)
v^{\eps}(0,\theta(\xi,\theta(x,t)))
S(\xi,\theta(x,t))\,d\xi\,
S(x,t),\vphi(x,t)>
$$
$$
=\int\limits_{I_+^{\eps_0}(i)}\int
\limits_{\om(\theta(x,t);0,t_{i-1}^*+\eps_{i-1}^+(\eps_0))}^{Q(\theta(x,t))}
b_s^{\eps}(\xi)\Bigl[c_r(\tau)v^{\eps}(0,\theta(\xi,\tau))
$$
$$
\times S(\xi,\tau)\Bigr]\Big|_{\tau=\theta(x,t)}S(x,t)\vphi(x,t)\,d\xi\,d(x,t)
$$
$$
=\sum\limits_{(j,q)\in J(i)}
\int\limits_{x_j^*-\eps_0}^{x_j^*+\eps_0}
\int\limits_{t_q^*-\eps_q^-(\eps_0)}^{t_q^*+\eps_q^+(\eps_0)}
\int\limits_{0}^{P(\xi,t)}
\Bigl[D(\xi,x,t)
\vphi(x,\tilde\om(x;0,\tilde
\om(\xi;0,t)))
$$
$$
-\chi_{[x_j-\eps,
x_j+\eps]\times[t_q^*-\eps_q^-(\eps),
t_q^*+\eps_q^+(\eps)]}(\xi,t)D(x_j,x,t_q^*)\vphi(x,\tilde\om(x;0,t_i^*))\Bigr]
b_s^{\eps}(\xi)v_q^{\eps}(0,t)
\,dx\,dt\,d\xi
$$
$$
+\sum\limits_{(j,q)\in J(i)}\int\limits_{x_j^*-\eps_0}^{x_j^*+\eps_0} b_s^{\eps}(\xi)\,d\xi
\int
\limits_{t_q^*-\eps_q^-(\eps_0)}^
{t_q^*+\eps_q^+(\eps_0)}v_q^{\eps}(0,t)
\,dt
 \int\limits_{0}^{P(\xi,t)}D(x_j,x,t_q^*)\vphi(x,
\tilde\om(x;0,t_i^*))\,dx,
$$
where $J(i)$ is the set of pairs $(j,q)$ such that $\om(x_j;0,t_q^*)=t_i^*$,
\begin{equation}
P(\xi,t)=
\cases{L
&if
$\theta(L,T)\ge\tilde\om(\xi;0,t)$,
 \cr \om(T;0,\tilde\om(\xi;0,t))
&if
$\theta(L,T)<\tilde\om(\xi;0,t)$, \cr}
\end{equation}
and
$$
D(\xi,x,t)
=
\biggl[ c_r(\theta(x,\tau))
S(\xi,\theta(x,\tau))
S(x,\tau)
$$
$$
\times\Bigl[(\d_t\theta)(\xi,\theta(x,
\tau))
(\d_t\theta)(x,\tau)\Bigr]^{-1}\biggr]
\bigg|_{\tau=\tilde\om(x;0,\tilde\om(\xi;0,t))}.
$$
It is immediate now that
$$
<v_i^{\eps}(x,t),\varphi(x,t)>
\longrightarrow_{\eps\to 0}
\sum\limits_{(j,q)\in J(i)}C_q
\int\limits_{0}^{P(x_j,t_i^*)} D(x_j,x,t_q^*)\vphi(x,
\tilde\om(x;0,t_i^*))\,dx.
$$
This  convergence is true due to
the second induction assumption, the condition (\ref{eq:eps_1}), the
continuity of $c_r,p,\la,$ and $\vphi$, and the uniform in $\eps$
boundedness of $\int_0^T|v^{\eps}(0,t)|\,dt$, that is proved
in Lemma 3.

Take a test function $\psi(t)\in\D(\R_+)$ and consider the
action $<v_i^{\eps}(0,t)$, $\psi(t)>$. Computation similar to the above
shows that
$$
v_i^{\eps}(0,t)\longrightarrow_{\eps\to 0}c_r(t_i^*)
\sum\limits_{(j,r)\in I(i)}\frac{S(x_j,t_i^*)}{(\d_t\om)(0;x_j,t_i^*)}
S(0,t_i^*)\de(t-t_i^*)
$$
$$
+c_r(t_i^*)\sum\limits_{(j,q)\in J(i)}C_q
\frac{S(x_j,t_i^*)}{(\d_t\om)(x_j;0,t_i^*)}S(0,t_i^*)
\de(t-t_i^*),
$$
thereby proving the required convergence.
\end{proof}

Lemma 4 is proved.

\section{Proof of Lemma 5}

We will need two facts from the theory of distributions.

\begin{fact}\label{thm:Horm}
(\cite{Horm}, 2.2.1)
If $u\in\D'(X)$ and every point of $X$ has a neighborhood on which
$u=0$, then $u\equiv 0$.
\end{fact}

\begin{fact}\label{thm:Joshi}
(\cite{Joshi}, 1.5.1)
If $(f_j)_{1\le j<\infty}\in L_{loc}^1(X)$ is a sequence which converges
almost everywhere to a function $f$, and there is a function
$g\in L_{loc}^1(X)$ such that $|f_j|\le g$ for all $j$, then
$f\in L_{loc}^1(X)$ and $f_j\to f$
in $\D'(X)$ as $j\to\infty$.
\end{fact}

Fact \ref{thm:Joshi} extends in an obvious way to  families of functions
$(f_{\eps})_{\eps>0}\in L_{loc}^1$.

To prove the lemma, it suffices to show
 that the conditions assumed in Fact~\ref{thm:Horm} are fulfilled.

Let us fix an arbitrary $(x,t)\in\Pi\setminus I_+$ and show that there
exists a neighborhood $Y\subset\Pi\setminus I_+$ of $(x,t)$ such that
the restriction of $v$ to $Y$ is equal to 0. We choose $Y$ such that
$\d Y\cap I_+=\emptyset$. This is possible because $\Pi\setminus I_+$
is open.
We now prove that $v=0$ on $Y$. By the definition of convergence in $\D'$,
if $v=\lim\limits_{\eps\to 0}v^{\eps}$ in $\D'(\Pi)$, then
$v=\lim\limits_{\eps\to 0}v^{\eps}$ in $\D'(Y)$. On the account of
Lemma 4, it suffices to prove the convergence of $v^{\eps}$ to 0
in $\D'(Y)$.

Let us check the conditions of Fact~\ref{thm:Joshi}
for $v^{\eps}$ on $Y$. The function $v^{\eps}$ is in $L_{loc}^1(Y)$
by Proposition 1. By Lemma 1, $v^{\eps}$ converges to 0 pointwise  on
$Y$ as $\eps\to 0$.
By the conditions imposed on $\la$, each component of $\Pi\setminus I_+$
is bounded. Since $Y$ is included in one of the components, $Y$ is bounded.
Clearly, there exists $T>0$ such that
$Y\subset\Pi^T$.
Since $\bigcap_{\eps>0}I_+^{\eps}=I_+$,
$\d Y\cap I_+=\emptyset$, and $\supp v^{\eps}\subset \overline{I_+^{\eps}}$
(see (\ref{eq:0})),
there exists $\tilde\eps$ such that
$v^{\eps}=0$ on $Y$ for all $\eps\ge\tilde\eps$.
 By (\ref{eq:1eps}), that was proved in Section 6, it follows
that $v^{\eps}=O\Bigl(\frac{1}{\eps}\Bigr)$ on $\Pi^T$. This implies the
uniform in $\eps$ estimate
$$
|v^{\eps}|\le\frac{C}{\tilde\eps}
$$
on $Y$.

Thus all conditions of
Fact~\ref{thm:Joshi}  hold for $v^{\eps}$ on $Y$.
Therefore $v^{\eps}\to 0$ in $\D'(Y)$ and
 $v=0$ on $Y$.

Since $(x,t)\in\Pi\setminus I_+$ is arbitrary, the lemma is true by
Fact~\ref{thm:Horm}.

\begin{rem}
Lemmas 4 and 5 show that $v^{\eps}$ are nets
converging to measures concentrated on $I_+$. From
 the  construction of $I_+$ (see Definition 2) it follows
that in general the density of
singularity curves increases as time progresses.
\end{rem}

\subsection*{Acknowledgments}
I am thankful to the members of the DIANA group  for their
hospitality during my stay at the Vienna
university.

\end{document}